\documentclass[times]{oupau}

\usepackage{indentfirst}  
\usepackage{mathrsfs}
\usepackage{mathtools}
\makeatletter
\renewenvironment{proof}[1][Proof]{%
  \par\medbreak
  \pushQED{\qed}%
  \noindent\textit{#1.}\hskip.5em\ignorespaces
}{%
  \popQED\medbreak
}
\makeatother 

\usepackage{algorithm}
\usepackage{algorithmic}
\usepackage{subcaption}
\usepackage{xcolor}
\usepackage{cite}
\usepackage{booktabs}
\usepackage{makecell}

\theoremstyle{remark}
\newtheorem{remark}{Remark}

\begin{document}

\runningheads{EDP-PI-DOSnet for Two-Phase Flows}{EDP-PI-DOSnet for Two-Phase Flows}

\title{EDP-PI-DOSnet: An Energy Dissipation Preserving Neural Solver for Two-Phase Flow Problems}

\author{Jizu Huang\affil{a,b},
Yue Qian\affil{a,b},  and  Tao Zhou\affil{a}}

\address{\affilnum{a}SKLMS, Academy of Mathematics and Systems Science, Chinese Academy of Sciences, Beijing, 100190, PR China.\\
\affilnum{b}School of Mathematical Sciences, University of Chinese Academy of Sciences, Beijing 100190, PR China.}

\begin{abstract}
Preserving the energy dissipation properties of gradient flow equations is important for the stability and reliability of numerical approximations, particularly in long-time simulations. In this paper, we propose a class of energy-dissipation-preserving physics-informed deep operator splitting networks (EDP-PI-DOSnets) for solving gradient flow equations. The proposed framework incorporates classical structure-preserving numerical schemes into neural network architectures, allowing the corresponding discrete energy dissipation laws to be preserved. In particular, suitable constraints are imposed on the convolution kernels to ensure that the learned discrete operators are symmetric and negative semidefinite. Furthermore, several network formulations are developed based on the Crank--Nicolson (CN), scalar auxiliary variable Crank--Nicolson (SAV/CN), and relaxed discrete variational derivative (R-DVD-1) schemes. Meanwhile, a physics-informed loss function based on the residuals of the governing equations is employed. Numerical experiments show that the predicted solutions are in good agreement with the reference solutions while satisfying the corresponding discrete energy dissipation laws. The results also suggest stable behavior during long-time recursive prediction and the potential applicability of the proposed framework to two-phase flow problems.
\end{abstract}

\keywords{Physics-informed learning; Gradient flows; Energy dissipation laws; Operator learning}


\maketitle


\section{Introduction}

Partial differential equations (PDEs) provide a fundamental mathematical framework for modeling a wide range of physical phenomena in science and engineering, including fluid dynamics, heat transfer, and multiphase flows \cite{chorin1967numerical,patankar1983calculation,drew1982mathematical}. Many PDEs possess intrinsic mathematical structures, such as conservation laws, energy dissipation laws, and other invariants. These structures not only reflect the underlying physical evolution of the corresponding systems but also provide important foundations for the construction and analysis of numerical methods. Classical numerical methods, such as the finite difference \cite{leveque2007finite}, the finite element \cite{zienkiewicz2013finite}, the finite volume \cite{leveque2002finite}, and spectral methods \cite{canuto2006spectral}, have been extensively developed to solve PDEs. In recent years, the approximation capabilities of neural networks have motivated the development of deep learning methods for PDEs \cite{sirignano2018dgm,raissi2019physics,yu2018deep,zang2020weak}. Among these approaches, physics-informed neural networks (PINNs) \cite{raissi2019physics} have received considerable attention. PINNs approximate the unknown solution using a neural network whose parameters are determined by minimizing a loss function incorporating the governing equations and the prescribed initial and boundary conditions. Despite their flexibility, neural-network-based methods do not generally guarantee the preservation of the intrinsic physical structures of the governing equations, such as conservation laws and energy dissipation properties. 
Therefore, developing neural-network-based numerical methods that achieve both high approximation accuracy and preservation of these intrinsic structures remains an important research objective.

Beyond solving individual PDE problems, operator learning seeks to approximate mappings between function spaces, enabling efficient predictions under varying initial conditions, boundary conditions, or physical parameters without retraining. Representative architectures include the deep operator network (DeepONet) \cite{lu2021learning}, the Fourier neural operator (FNO) \cite{li2020fourier}, and the deep operator-splitting network (DOSnet) \cite{lan2023dosnet}. Existing operator-learning methods can be broadly categorized into data-driven and physics-informed approaches. Data-driven methods learn solution operators from paired input-output data \cite{li2020neural,lan2023dosnet,tripura2023wavelet,chen2024positional,hao2023gnot,ye2024pdeformer,wu2024transolver}, whereas physics-informed approaches incorporate the governing equations and associated constraints into the training process, thereby reducing the reliance on labeled solution data \cite{wang2021learning,li2024physics,gao2021phygeonet,navaneeth2024physics}. In particular, PI-DOSnet \cite{huang2026pidosnet} integrates the operator-splitting architecture of DOSnet \cite{lan2023dosnet} with physics-informed training to approximate the solution operators of nonlinear evolution PDEs over continuous time intervals. Nevertheless, neither data-driven approximation nor physics-informed approaches inherently guarantees the preservation of the  underlying physical structures of the governing equations. In particular, ensuring that learned solution operators satisfy energy dissipation laws remains an important challenge for the reliable long-time prediction of dissipative systems.

To preserve the underlying physical structures  of PDEs, based on traditional discrete schemes, various structure-preserving numerical methods have been developed. These methods aim to retain the corresponding properties of the governing equations at the discrete level, thereby improving the reliability of numerical simulations, particularly for long-time evolution. Representative approaches include symplectic and multisymplectic integrators for Hamiltonian systems \cite{sanz1992symplectic,bridges2001multi}, compatible and mimetic discretizations that preserve structural relations among differential operators \cite{arnold2006finite,hyman2002mimetic}, and invariant-preserving schemes for conserved quantities \cite{dahlby2011general}. For dissipative systems, particularly gradient flows, considerable attention has been devoted to constructing numerical schemes that preserve the underlying energy dissipation law. Convex-splitting methods achieve unconditional energy stability by treating the convex and concave parts of the energy implicitly and explicitly, respectively \cite{elliott1993global,eyre1998unconditionally}. Stabilized semi-implicit methods introduce suitable stabilization terms to control nonlinear contributions and establish energy stability under appropriate conditions \cite{shen2010numerical}. Discrete variational derivative (DVD) methods enforce a discrete chain rule, leading to schemes that satisfy a discrete dissipation law for the original energy \cite{du1991numerical,huang2025high}. Alternatively, invariant energy quadratization (IEQ) and scalar auxiliary variable (SAV) methods introduce auxiliary variables to facilitate the construction of linear and often decoupled schemes that are unconditionally energy stable with respect to a modified energy \cite{yang2016linear,shen2018scalar}. These energy-stable approaches have also been extended to coupled phase-field and fluid models of two-phase flows \cite{li2020sav,yang2021novel}.

Building on classical structure-preserving numerical methods, recent studies have explored incorporating conservation and dissipation properties into neural models of dynamical systems. For Hamiltonian systems, Hamiltonian neural networks learn Hamiltonian functions to represent the underlying dynamics, while symplectic networks construct evolution maps that preserve the underlying symplectic structure \cite{greydanus2019hamiltonian,chen2019symplectic,jin2020sympnets}. For dynamical systems possessing first integrals, invariant-preserving neural networks employ structured parameterizations or projection techniques to maintain the corresponding conserved quantities along predicted trajectories \cite{matsubara2022finde,cardoso2025exactly}. For dissipative systems, neural approaches to promoting or preserving energy dissipation can be broadly divided into training-based constraints and structure-preserving constructions. The former incorporate energy penalties or physically motivated regularization into training objectives to improve the physical consistency of predictions, but generally do not guarantee discrete energy dissipation \cite{kutuk2025energy,li2026structure}. In contrast, the latter embed structure-preserving mechanisms directly into neural architectures or numerical time discretizations. For instance, EStable-Net incorporates auxiliary-variable-based network blocks to guarantee the dissipation of a constructed discrete energy, although it relies on supervised learning \cite{wu2023energy}. Other approaches combine neural representations with discrete-gradient methods or energetic variational formulations \cite{matsubara2020deep,hu2024energetic}. 
More recently, neural approaches based on SAV, such as EDE-DeepONet \cite{zhang2024energy} and the method of Luo et al. \cite{luo2026sav}, evolve neural parameters or coefficients through successive time integration, which achieve modified energy dissipation.
Despite these advances, developing physics-informed neural operators that simultaneously guarantee discrete energy dissipation, accommodate multiple initial conditions, and support reliable long-time prediction without retraining remains an important challenge.

Motivated by this challenge, we build upon PI-DOSnet \cite{huang2026pidosnet} to develop a physics-informed operator-learning framework with guaranteed discrete energy dissipation. PI-DOSnet employs physics-informed training based on the governing equations to learn reusable solution operators that accommodate multiple initial conditions and enable long-time prediction. However, its original operator-splitting architecture, which alternates between learned linear evolution and nonlinear subproblems, does not inherently guarantee discrete energy dissipation for dissipative systems. To overcome this limitation, we impose structural constraints on the learned linear operator to ensure its negative semidefiniteness, thereby making the linear evolution compatible with energy dissipation. Furthermore, we exploit the blockwise architecture of PI-DOSnet to incorporate classical energy-stable semi-implicit schemes, specifically the SAV and DVD methods, into the network evolution. These developments lead to the proposed Energy-Dissipation-Preserving PI-DOSnet (EDP-PI-DOSnet) framework, which rigorously satisfies a discrete energy dissipation law while retaining the operator-learning capabilities of PI-DOSnet. We conduct numerical experiments on representative gradient-flow and two-phase-flow problems to evaluate the accuracy, energy stability, and robustness of the proposed framework. The results demonstrate that EDP-PI-DOSnet effectively preserves energy dissipation and achieves more accurate and reliable long-time predictions than the original PI-DOSnet.

The remainder of this paper is organized as follows. Section \ref{sec:pre} reviews the PI-DOSnet framework, summarizes an existing energy-constrained extension of PI-DOSnet, and introduces the classical energy-stable formulations employed in this work. Section \ref{sec:edp-pidos} presents the proposed EDP-PI-DOSnet for solving gradient-flow problems. Section \ref{sec:experiments} provides numerical experiments to demonstrate the effectiveness and robustness of the proposed method. Finally, Section \ref{sec:conclusion} concludes the paper.

\section{Premiliaries} \label{sec:pre}

\subsection{Review of the PI-DOSnet}\label{subsec:pidos}

The PI-DOSnet \cite{huang2026pidosnet} is designed to solve a general autonomous evolution PDE of the form:
\begin{equation}\label{eq:original_equation}
\begin{aligned}
&u_t = \mathcal{L}u+\mathcal{N}u, \quad &\text{in} \ (0,T]\times \Omega,\\
&u\vert_{t=0}=u_0,&\text{in} \ \Omega,\\
		&\mathcal{B}u=0,&\text{on} \ \partial\Omega\times(0,t^\star],
\end{aligned}
\end{equation}
where $\mathcal B$ is the boundary operator, $u_0$ denotes the initial condition, $u:=u(t,\bm x)$ is the unknown solution,  $\mathcal{L}$, and $\mathcal{N}$ are the linear and nonlinear differential operators, respectively.
The operators $\mathcal{B}$, $\mathcal{L}$, and $\mathcal{N}$ are assumed to be independent of the time variable $t$. The goal of PI-DOSnet is to learn the solution operator
\begin{equation}
	\mathcal{G}:\mathcal{U}\rightarrow \mathcal{S}\subset\mathcal{X},
\end{equation}
which maps a given initial condition to its corresponding solution trajectory. More specifically, for any initial condition \(u_0\in\mathcal{U}\) and any time \(t\in[0,T]\), the network output provides an approximation of the exact solution.

The architecture of PI-DOSnet is inspired by the operator-splitting framework of DOSnet \cite{lan2023dosnet}, but is reformulated within a physics-informed manner. As in DOSnet, the evolution operator is decomposed into linear and nonlinear components. The linear component is represented by learnable convolutional operators, whereas the nonlinear component is modeled through the solution operator of the corresponding nonlinear subproblem. 
The primary distinction between PI-DOSnet and DOSnet lies in the treatment of the linear evolution operator. Rather than directly learning the exponential operator \(e^{dt\mathcal{L}}\), PI-DOSnet approximates it using a low-order Taylor expansion,
\begin{equation}\label{eq:taylor_expansion_pidosnet}
e^{dt\mathcal{L}}
\approx
\mathcal{I}+dt\mathcal{L}_{\bm\theta}
+\frac{dt^2}{2}\mathcal{L}_{\bm\theta}^2,
\end{equation}
where $dt$ is a parameter corresponding to the time step size and $\mathcal{L}_{\bm\theta}$ denotes a  learned convolutional operator. 
This approximation makes the learned linear operator appear explicitly in the PDE residual, thereby enabling the network to be trained using physical constraints alone, without requiring reference solution data.

For the nonlinear component, PI-DOSnet adopts the same philosophy as DOSnet \cite{lan2023dosnet}. Rather than employing a generic activation function, the nonlinear layer is constructed from the analytical solution operator associated with the split nonlinear subproblem
\begin{equation}
u_t=\mathcal{N}u.
\end{equation}
Consequently, the nonlinear layer is problem-dependent and determined by the specific form of \(\mathcal{N}\). It may be viewed as a dynamics-informed activation function that incorporates the intrinsic nonlinear evolution of the underlying PDE. As a result, the applicability of PI-DOSnet relies on the availability of an explicit or analytically tractable solution operator for the nonlinear subproblem.

\begin{figure}[!htbp]
	\centering
	\includegraphics[scale=0.28]{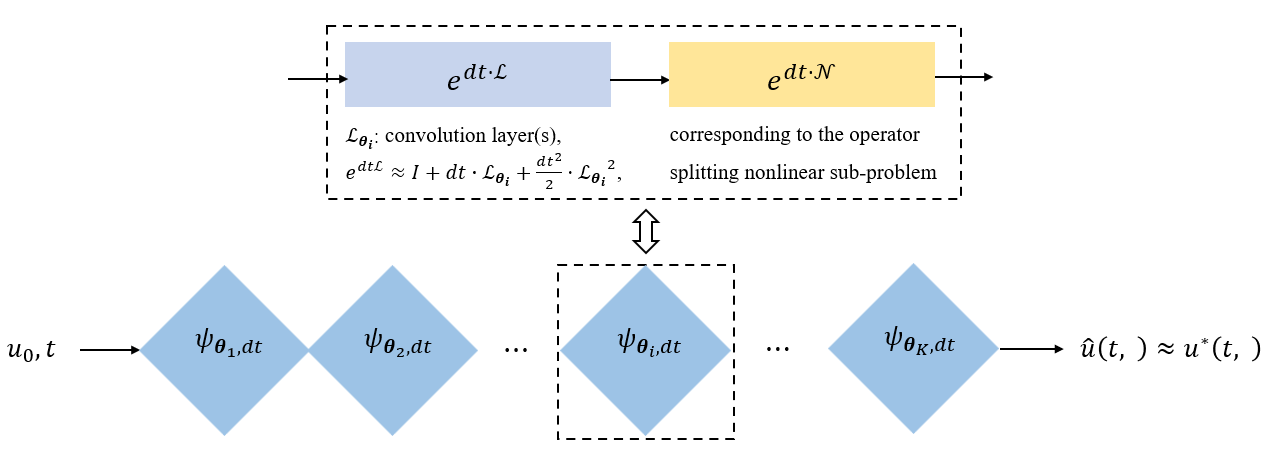}
	\caption{Architecture of PI-DOSnet.}
	\label{fig:pidosnet_blocks}
\end{figure}

As illustrated in Figure~\ref{fig:pidosnet_blocks},
PI-DOSnet consists of \(K\) splitting blocks, each comprising a linear layer followed by a nonlinear layer. Given a target time $t$, we define the parameter in \eqref{eq:taylor_expansion_pidosnet} as
\begin{equation}
	dt=\frac{t}{K}.
\end{equation}
The $i$th splitting block of PI-DOSnet is then defined by
\begin{equation}\label{PI-DOSnet-ith-block}
	\hat{u}(t_i,\bm x)=\psi_{\bm{\theta}_i,dt}(\hat{u}(t_{i-1},\bm x))=e^{dt\mathcal{N}}\circ\left(\mathcal{I}+dt\mathcal{L}_{\bm{\theta}_i}+\frac{dt^2}{2}\mathcal{L}_{\bm{\theta}_i}^2\right)\left(\hat{u}(t_{i-1},\bm x)\right),
\end{equation}
where $t_i=i \times dt$ and $\mathcal{L}_{\bm{\theta}_i}$ denotes the learned convolutional operator. By successively composing all splitting blocks, the PI-DOSnet approximation can be expressed as
\begin{equation}\label{eq:pidosnet_mapping}
	\hat{u}(t,\bm{x})=\Phi_{\bm{\theta}}(u_0,t):=\psi_{\bm{\theta}_K,dt}\circ\psi_{\bm{\theta}_{K-1},dt}\circ\cdots\circ\psi_{\bm{\theta}_1,dt}(u_0).
\end{equation}
This recursive composition approximates the solution trajectory by repeatedly applying the same operator-splitting principle, thereby enabling the network to capture the time evolution of the underlying dynamical system over the interval \([0,t]\).

Let $\{u_0^{(n)}\}_{n=1}^{N_u}$ denote the training set of initial conditions. For each $u_0^{(n)}$, we sample a set of interior collocation points $\{(t_{r,j}^{(n)},\bm{x}_{r,j}^{(n)})\}_{j=1}^{N_r}$ and boundary collocation points $\{(t_{b,j}^{(n)},\bm{x}_{b,j}^{(n)})\}_{j=1}^{N_b}$. 
The spatial collocation points are selected from a fixed uniform grid $\{\bm{x}_i\}_{i=1}^{N_x}$, whereas the temporal collocation points are randomly sampled from the training interval $[0,T]$.
The loss function of PI-DOSnet is defined as
\begin{equation}\label{eq:pidosnet_loss}
	\mathfrak{L}(\bm{\theta})=\lambda_r\mathfrak{L}_r+\lambda_b\mathfrak{L}_b,
\end{equation}
where \(\mathfrak{L}_r\) and \(\mathfrak{L}_b\) denote the PDE residual loss and boundary-condition loss, respectively:
\begin{equation}
	\mathfrak{L}_r=\frac{1}{N_uN_r}\sum_{n=1}^{N_u}\sum_{j=1}^{N_r}\left|\left(\hat{u}_{r,j}^{(n)}\right)_t-\mathcal{L}\hat{u}_{r,j}^{(n)}-\mathcal{N}\hat{u}_{r,j}^{(n)}\right|^2,
\end{equation}
and
\begin{equation}
	\mathfrak{L}_b=\frac{1}{N_uN_b}\sum_{n=1}^{N_u}\sum_{j=1}^{N_b}\left|\mathcal{B}\hat{u}_{b,j}^{(n)}\right|^2.
\end{equation}
Here, $\lambda_r$ and $\lambda_b$ are prescribed weighting parameters that balance the contributions of the PDE residual and boundary-condition losses. 
In the evaluation of the loss function, the time derivative $u_t$ is computed using automatic differentiation, while the spatial derivatives of $u$ are approximated by finite difference schemes.
The trainable parameters are then obtained by minimizing the total loss $\mathfrak{L}(\bm{\theta})$. 

Although PI-DOSnet is a physics-informed network, it does not automatically preserve the implicit physical constraints inherent in the governing PDEs. We take the Allen--Cahn (AC) equation as an example to illustrate this point. The AC equation is defined as 
\begin{equation}\label{eq:AC}
    	\frac{\partial u}{\partial t} = - \frac{\delta \mathcal{E}(u)}{\delta u}:=\epsilon^2 \Delta u - F'(u), 
    	\qquad 
    	\mathcal{E}(u) = \frac{\epsilon^2}{2} \langle \nabla u, \nabla u \rangle + \mathcal{E}_1(u),
    \end{equation}
    where $\epsilon$ denotes the interface width, $F(u)=(u^2-1)^2$ is the double well potential, $\langle \phi, \psi \rangle = \int_\Omega \phi \psi \, \mathrm{d}\bm x$ denotes the standard $L^2$ inner product, and $\mathcal{E}_1(u) = \int_\Omega F(u) \, \mathrm{d}\bm x$ represents the bulk energy. Under periodic boundary conditions, the energy dissipation law associated with the AC equation takes the form 
    \begin{equation}
	\frac{\mathrm{d}\mathcal{E}}{\mathrm{d}t} = - \|u_t\|_{L^2}^2,
\end{equation}
where $\|u_t\|_{L^2}:=\langle u_t, u_t \rangle^{1/2}$
As observed in the numerical results of \cite{huang2026pidosnet}, predictions produced by PI-DOSnet may exhibit energy growth between consecutive time steps, thereby violating the underlying energy dissipation law. This raises a key issue: how can the architecture of PI-DOSnet be enhanced so that the implicit physical constraints inherent in the governing PDEs are properly respected? Another motivation for preserving such implicit physical constraints stems from the long-time prediction performance of PI-DOSnet. As demonstrated by the numerical results in \cite{huang2026pidosnet}, when solving the AC equation, a PI-DOSnet trained on the time interval $[0,T]$ and subsequently applied to predict the solution over $[0,10T]$ typically requires two or three rounds of retraining to achieve the desired accuracy. Such frequent retraining significantly increases the computational cost and limits the efficiency of PI-DOSnet for long-time simulations. This observation further motivates the development of architectures that preserve implicit physical constraints in order to improve long-time prediction stability. 
In the remainder of this paper, we use the AC equation as a representative example to illustrate our proposed improvements based on penalty formulations and traditional structure-preserving numerical methods.

\subsection{PI-DOSnet with Penalty Function}\label{subsec:penalty_pidos}
When incorporating physical constraints into neural network models, such as the energy dissipation law associated with gradient flow problems, a natural and widely adopted strategy is the penalty method. The central idea is to augment the loss function with additional terms that penalize deviations from the prescribed physical behavior. Following this philosophy, while solving gradient flow problems, i.e. the AC equation \eqref{eq:AC}, EStable-Net \cite{wu2023energy} explicitly incorporates a discrete energy-dissipation penalty into its training objective:
\begin{equation}
	\mathcal{L}_{\Delta\mathcal{E}}(\bm\theta) = \beta \sum_{i=1}^{K} \big[ \mathcal{E}(u^{i}) - \mathcal{E}(u^{i-1}) - \delta \mathcal{E}^{i-1}(u^{i-1}, u^{i}) \big]^2,
\end{equation}
where $u^{i}\approx u(t_i,\bm x)$ and $\mathcal{E}(u^{i}) - \mathcal{E}(u^{i-1})$ represents the actual energy change predicted by the network between consecutive states, while $\delta\mathcal{E}^{i-1}$ denotes the theoretical energy variation dictated by the underlying physical model. This penalty term serves two purposes: it discourages unphysical energy growth and simultaneously guides the learned energy decay rate toward the theoretically prescribed one. Consequently, the discrete energy variation can be viewed as a numerical approximation of the continuous energy dissipation law.

For the Allen--Cahn (AC) equation, approximating the time derivative by a finite difference over a time step \( dt\) yields the corresponding discrete energy variation
\begin{equation}
    	\delta\mathcal{E}^{i-1}(u^{i-1}, u^{i}) = - \frac{1}{dt} \left\| u^{i} - u^{i-1} \right\|_{L^2}^2.
\end{equation}
By embedding this discrete counterpart of the continuous energy law into the optimization objective, EStable-Net  effectively integrates physical principles into the learning process. As a result, the resulting neural network not only achieves accurate approximations of the solution but also preserves the physically meaningful energy-dissipation behavior, thereby enhancing both robustness and physical consistency.

To enforce physical constraints, it is straightforward to augment PI-DOSnet with a penalty term that encodes the desired physical law. For the AC equation, the resulting loss function takes the form
\begin{equation}\label{eq4:energy_penalty_loss}
\begin{aligned}
    & \mathfrak{L}(\bm{\theta}) = \lambda_r \mathfrak{L}_r + \lambda_b \mathfrak{L}_b + \lambda_e \mathfrak{L}_{\mathcal{E}_t}, \\
    & \mathfrak{L}_{\mathcal{E}_t} = \frac{1}{N_u} \frac{1}{N_t} \sum_{n=1}^{N_u} \sum_{i=1}^{N_t} \left| \frac{\mathrm{d} \mathcal{E}\left(\hat{ u}^{i,(n)}\right)}{\mathrm{d} t} + \left\| \left(\hat{ u}_t\right)^{i,(n)} \right\|_{L^2}^2 \right|^2,
\end{aligned}
\end{equation}
where $\hat{ u} = \Phi_{\bm{\theta}_T}(u_0, t)$ denotes the PI-DOSnet prediction, representing the solution on the spatial grid at time $t$. Since all time derivatives appearing in $\mathfrak{L}_{\mathcal{E}_t}$ can be computed efficiently via automatic differentiation, the energy-dissipation constraint can be incorporated into the training process in a fully differentiable manner without requiring additional data. Consequently, PI-DOSnet can be trained to satisfy the underlying energy law while maintaining the data-free and physics-informed nature of the framework.

However, it should be noted that the penalty-based approach has inherent limitations. The effectiveness of the constraint critically depends on the weight $\lambda_e$ and the adequacy of the sampling in both space and time. If the weight is too small, the constraint on energy stability may not be enforced adequately; conversely, excessive weight can dominate residual loss and impede the ability of the network to accurately approximate the PDE solution. Furthermore, this approach does not guarantee exact enforcement of energy laws at every time step or at every spatial collocation point; it only biases the optimization towards energy-consistent trajectories at the sampled temporal and spatial points. Therefore, while the penalty method offers a practical mechanism to incorporate energy considerations, careful tuning and validation are necessary to maintain a proper balance between fidelity to the governing PDE and adherence to the underlying energy principle across the sampled configurations. 

\subsection{Introduction to Classical Energy-Stable Methods}

    For gradient flow problems, there exist many well-established classical energy-stable numerical methods that preserve the energy dissipation law. We take the AC equation as an example to briefly introduce several representative approaches, whose underlying ideas will be employed in the next section to develop a penalty-free operator learning solver that preserves the energy dissipation law.
     \begin{itemize}
    	\item \textbf{Convex Splitting (CS) Methods} \cite{elliott1993global, eyre1998unconditionally}: These methods split the energy functional into the difference of two convex components,
    	\begin{equation*}
    		\mathcal{E}(u) = \mathcal{E}_c(u) - \mathcal{E}_e(u),
    	\end{equation*}
    	and treat $\mathcal{E}_c(u)$ implicitly and $\mathcal{E}_e(u)$ explicitly in the discretization. A first-order discrete scheme reads
    	\begin{equation*}
    		\frac{u^{n+1} - u^n}{\mathrm{d}t} = -\nabla \mathcal{E}_c(u^{n+1}) + \nabla \mathcal{E}_e(u^n),
    	\end{equation*}
    	which is unconditionally energy-stable, i.e., the discrete energy decreases for any $\Delta t>0$. However, such schemes typically require iterative nonlinear solvers at each time step.
    
    	\item \textbf{Stabilized Methods} \cite{shen2010numerical}: These approaches introduce an artificial stabilization term to balance the explicit treatment of nonlinear terms. For the AC equation, a classical stabilized scheme is
    	\begin{equation*}
    		\frac{u^{n+1}-u^n}{\mathrm{d}t} = \epsilon^2 \Delta u^{n+1} - F'(u^n) - S (u^{n+1}-u^n),
    	\end{equation*}
    	where $S$ is the stabilization constant. If $\|F''(u)\|_\infty \le a$, unconditional stability is guaranteed for $S \ge a/2$. The main advantage is that only a linear system needs to be solved at each time step, leading to improved computational efficiency.
    
    	\item \textbf{Scalar Auxiliary Variable (SAV) Methods} \cite{shen2018scalar}: SAV refines IEQ by introducing a scalar auxiliary variable
    	\begin{equation*}
    		r(t) = \sqrt{\int_\Omega F(u) \, \mathrm{d}\bm x + C_0},
    	\end{equation*}
    	to overcome variable-coefficient difficulties. The AC equation becomes
    	\begin{equation}\label{eq:sav}
    		\begin{aligned}
    			& u_t = \epsilon^2 \Delta u - \frac{r(t)}{\sqrt{\int_\Omega F(u) \, \mathrm{d}\bm x + C_0}} F'(u), \\
    			& r_t = \frac{1}{2 \sqrt{\int_\Omega F(u) \, \mathrm{d}\bm x + C_0}} \int_\Omega F'(u) u_t \, \mathrm{d}\bm x.
    		\end{aligned}
    	\end{equation}
        The total energy is modified as
        \begin{equation*}
    		\mathcal{E}_{SAV}(u) = \frac{\epsilon^2}{2} \langle \nabla u, \nabla u \rangle + r^2.
    	\end{equation*}
    	Combining SAV with a Crank--Nicolson scheme yields a second-order semi-implicit time discretization (SAV/CN):
    	\begin{subequations}\label{eq:cnsav}
    		\begin{align}
    			& \frac{u^{n+1} - u^n}{\mathrm{d}t} = \epsilon^2 \Delta u^{n+1/2} - r^{n+1/2} \bar{b}^{n+1/2}, \\
    			& r^{n+1} - r^n = \frac{1}{2} \langle \bar{b}^{n+1/2}, u^{n+1} - u^n \rangle,
    		\end{align}
    	\end{subequations}
    	with $\bar{b}^{n+1/2}$ computed from a second-order explicit predictor. This scheme guarantees unconditional modified energy dissipation. Although implicit, the SAV/CN system can be decoupled into two linear systems, solved efficiently with standard solvers.
    
    	\item \textbf{Relaxed Discrete Variational Derivative (R-DVD) Methods} \cite{huang2025high}: DVD methods provide a rigorous framework for enforcing discrete energy dissipation by constructing a discrete variational derivative $\mu$ such that
    	\begin{equation*}
    		\mathcal{E}^{n+1} - \mathcal{E}^n = \langle u^{n+1} - u^n, \mu[u^{n+1}, u^n] \rangle.
    	\end{equation*}
    	Inspired by SAV, R-DVD relaxes strict energy preservation by introducing a modified energy $\hat{\mathcal{E}}$ and a scalar auxiliary variable $r = (\mathcal{E}_1(u) + C_0)^{1/m}$, yielding
    	\begin{equation*}
    		\hat{\mathcal{E}} = \frac{\epsilon^2}{2} \langle \nabla u, \nabla u \rangle + r^m,
    	\end{equation*}
    	and a semi-implicit scheme
    	\begin{equation}\label{eq:rdvd_scheme}
    		\frac{u^{n+1} - u^n}{\mathrm{d}t} = - \mu[u^{n+1}, u^n], \quad r^{n+1} = r^n + \frac{\langle F'(\bar{u}^{n+1/2}), u^{n+1} - u^n \rangle}{m r^{m-1}(\bar{u}^{n+1/2})},
    	\end{equation}
    	with $\bar{u}^{n+1/2}$ a second-order explicit approximation. This approach preserves unconditional energy stability and transforms the original nonlinear system into a computationally efficient linear system. Notably, for $m=2$, R-DVD reduces to classical SAV, while for $m=1$, it yields the fully decoupled R-DVD-1 scheme, where each time step requires solving only a single linear system for $u$, with explicit updates for $r$, combining high efficiency and strict energy stability.
    \end{itemize}

\section{Energy Dissipation Preserving PI-DOSnet}\label{sec:edp-pidos}
For nonlinear gradient flow problems, such as the AC equation, the ability to preserve the energy dissipation property (or energy stability) is a fundamental criterion for assessing the long-time performance of a numerical method. As discussed in the previous section, many classical numerical schemes, including the CS, SAV, and R-DVD methods, generate numerical solutions that satisfy a discrete energy dissipation law.
In contrast, PI-DOSnet generally does not preserve this property for two main reasons. First, the linear operator learned by PI-DOSnet is not guaranteed to be symmetric negative semidefinite. Although the training process enables PI-DOSnet to accurately approximate the underlying linear operator, such an approximation does not necessarily retain its structural properties. In particular, for the AC equation, the numerical results reported in \cite{huang2026pidosnet} show that the eigenvalues of the learned operator are not confined to the real axis, suggesting that the learned operator may fail to satisfy the negative-semidefiniteness required for energy dissipation. This deficiency can significantly compromise the energy stability of PI-DOSnet.
Second, the simple linear-nonlinear splitting strategy employed in PI-DOSnet does not, by itself, guarantee energy dissipation. As demonstrated by classical energy-stable schemes such as the SAV method, preserving the energy dissipation property generally requires a carefully designed discretization of both the linear and nonlinear terms. Consequently, the operator-splitting mechanism adopted in the original PI-DOSnet framework is insufficient to ensure a discrete energy dissipation law.

To develop an energy dissipation preserving PI-DOSnet, we address the two fundamental issues identified above. On the one hand, we impose appropriate constraints on the kernels of the convolution operators so that the resulting operators are symmetric negative semidefinite. We establish the corresponding theoretical framework, analyze the properties of the constrained convolution operators, and rigorously prove their symmetric negative-semidefinite structure.
On the other hand, we improve the operator-splitting mechanism of PI-DOSnet by incorporating classical semi-implicit time discretization schemes, such as the SAV and R-DVD methods. By embedding these structure-preserving numerical schemes into the network architecture, we construct PI-DOSnets that inherit the energy dissipation properties of the underlying numerical methods. Finally, we demonstrate that the resulting neural solver achieves energy stability by construction, eliminating the need for additional energy-penalty terms in the loss function.

\subsection{Negative-Semidefinite Convolutional Modules}

\subsubsection{Single-Layer Convolutional Operator}
A single-layer convolutional operator, or equivalently a convolutional layer, can be interpreted as a linear operator acting on the vectorized input tensor. To illustrate this relationship, we consider a one-dimensional convolutional layer with one input channel and one output channel. Let the kernel size be $2a+1$, and denote the convolution kernel by
\begin{equation}\label{eq:centered_kernel}
	\bm{w}=[w_{-a},w_{-a+1},\ldots,w_0,\ldots,w_{a-1},w_a]^T\in\mathbb{R}^{2a+1}.
\end{equation}
Assume that circular padding is used. For an input vector
$
	\bm{y}=[y_1,y_2,\ldots,y_n]^T\in\mathbb{R}^n,$ 
the output of the convolutional operator,  denoted as $\bm{z}=\mathcal{L}_w\bm{y}$, is given in the order of components by
\begin{equation}
	z_i=\sum_{r=-a}^{a}w_r y_{i+r},\qquad i=1,2,\ldots,n,
\end{equation}
where the index $i+r$ is understood periodically, i.e., $y_{i+r}=y_{i+r+n}$ or $y_{i+r}=y_{i+r-n}$ whenever the index falls outside $\{1,2,\ldots,n\}$. Therefore, the convolutional operation can be equivalently represented by a circulant matrix $\bm{L}\in\mathbb{R}^{n\times n}$:
\begin{equation}
	\mathcal{L}_w\bm{y}=\bm{L}\bm{y},
\end{equation}
where $\bm{L}$ is defined by 
\begin{equation} L_{ij}=\sum_{r=-a}^{a}w_r\,\mathbf{1}_{\{j\equiv i+r\;(\mathrm{mod}\; n)\}},\qquad i,j=1,2,\ldots,n. 
\end{equation} 
Here, $\mathbf{1}_{\{\cdot\}}$ denotes the indicator function. In other words, the entry $L_{ij}$ is equal to the kernel coefficient $w_r$ whenever the $j$th input component contributes to the $i$th output component through the offset $r$ under periodic indexing, and it is zero otherwise. For a more intuitive view of this matrix representation, we further display the circulant structure of $\bm{L}$ in a schematic form, namely,
\begin{equation}
	\bm{L}=
	\begin{bmatrix}
		w_0 & w_1 & \cdots & w_a & 0 & \cdots & 0 & w_{-a} & \cdots & w_{-1} \\
		w_{-1} & w_0 & \cdots & w_{a-1} & w_a & \cdots & 0 & 0 & \cdots & w_{-2} \\
		w_{-2} & w_{-1} & \cdots & w_{a-2} & w_{a-1} & \cdots & 0 & 0 & \cdots & w_{-3} \\
		\vdots & \vdots & \ddots & \vdots & \vdots & \ddots & \vdots & \vdots & \ddots & \vdots \\
		w_1 & w_2 & \cdots & 0 & 0 & \cdots & 0 & 0 & \cdots & w_0
	\end{bmatrix}.
\end{equation}
Here, the displayed matrix is intended to illustrate the circulant pattern induced by circular padding, while the precise position of each kernel coefficient is determined by the entrywise definition of $L_{ij}$ given above.

Thus, although the convolutional layer is implemented as an operation on tensor-shaped data, after vectorization it is equivalent to multiplying the input vector by a structured matrix. In the case of circular padding, this matrix is circulant. In what follows, we use $\mathcal{L}_w$ to denote the convolutional operator and $\bm{L}$ to denote its equivalent matrix representation. The notation $\mathcal{L}\bm{y}$ is therefore understood in this formal sense, with the reshaping between tensor and vector forms omitted for simplicity.
A convolutional layer $\mathcal{L}_w$ with a negative-semidefinite structure means that its equivalent matrix representation $\bm{L}$ is symmetric negative semidefinite. Starting from a general convolutional layer, a natural way to impose this property is to restrict the convolution kernel parameters. 

To ensure that the equivalent circulant matrix $\bm{L}$ is negative semidefinite, we explicitly impose the following constraints on the kernel coefficients:
\begin{equation}\label{nsd_conditions}
	\left\{
	\begin{aligned}
		& w_i=w_{-i},\qquad i=1,2,\ldots,a,\\
		& w_0=-2\sum_{i=1}^{a}|w_i|.
	\end{aligned}
	\right.
\end{equation}
The first condition enforces the symmetry of the convolution kernel and therefore ensures that the associated circulant matrix $\bm{L}$ is  symmetric. The second condition controls the diagonal coefficient relative to the off-diagonal weights. In particular, since each row of $\bm{L}$ contains the diagonal entry $w_0$ together with the off-diagonal coefficients generated by $w_{\pm i}$, $i=1,2,\ldots,a$, these constraints are designed to place the spectrum of $\bm{L}$ on the non-positive real axis. To make this mechanism more transparent, we recall the Gershgorin circle theorem and a standard corollary linking (strict) diagonal dominance of a symmetric matrix to the non-positivity of its eigenvalues.

\begin{theorem}[Gershgorin Circle Theorem \cite{horn2012matrix}]
\label{Gershgorin Circle Theorem}
	Let $\bm{A}=(a_{ij})$ be an $n\times n$ real matrix. If $\lambda$ is an eigenvalue of $\bm{A}$, then there exists an index $i\in\{1,2,\ldots,n\}$ such that
	\begin{equation}
		|\lambda-a_{ii}|\leq \sum_{j\neq i}|a_{ij}|.
	\end{equation}
	Equivalently, the spectrum of $\bm{A}$, denoted by $\sigma(\bm{A})$, is contained in the union of the Gershgorin disks:
	\begin{equation}
		\sigma(\bm{A})\subseteq \bigcup_{i=1}^{n}D_i,
	\end{equation}
	where
	\begin{equation}
		D_i=\left\{z\in\mathbb{C}: |z-a_{ii}|\leq \sum_{j\neq i}|a_{ij}|\right\}
	\end{equation}
	is the $i$th Gershgorin disk.
\end{theorem}

\begin{corollary}\label{coro_nsd}
	Let $\bm{A}=(a_{ij})$ be an $n\times n$ real symmetric matrix. Suppose that $\bm{A}$ is weakly diagonally dominant, namely,
	\begin{equation}
		|a_{ii}|\geq \sum_{j\neq i}|a_{ij}|,\qquad i=1,2,\ldots,n.
	\end{equation}
	If all diagonal entries are non-positive, i.e.,
	\begin{equation}
		a_{ii}\leq 0,\qquad i=1,2,\ldots,n,
	\end{equation}
	then $\bm{A}$ is negative semidefinite.
\end{corollary}

Corollary~\ref{coro_nsd} follows directly from Theorem~\ref{Gershgorin Circle Theorem}.
Under the constraints \eqref{nsd_conditions}, the equivalent matrix $\bm{L}$ satisfies all the assumptions of Corollary \ref{coro_nsd}. Indeed, the condition $w_i=w_{-i}$ makes $\bm{L}$ symmetric, while $w_0=-2\sum_{i=1}^{a}|w_i|$ gives $w_0\leq 0$. In addition, for each row of the circulant matrix $\bm{L}$, the sum of the absolute values of the off-diagonal entries is
\begin{equation}
	\sum_{j\neq i}|L_{ij}|=2\sum_{k=1}^{a}|w_k|=|w_0|.
\end{equation}
Hence, $\bm{L}$ is weakly diagonally dominant with non-positive diagonal entries. Therefore, by Corollary \ref{coro_nsd}, $\bm{L}$ and ${\mathcal L}_w$ are negative semidefinite.

Similarly, a positive-semidefinite convolutional layer can be constructed by reversing the sign of the kernel coefficient while keeping the symmetry constraint unchanged. More precisely, under the same indexing convention, we impose
\begin{equation}\label{psd_conditions}
	\left\{
	\begin{aligned}
		& w_i=w_{-i},\qquad i=1,2,\ldots,a,\\
		& w_0=2\sum_{i=1}^{a}|w_i|.
	\end{aligned}
	\right.
\end{equation}

\subsubsection{Multi-Layer Convolutional Operator}
We first consider an $m$-layer single-channel convolutional module. Let $\bm{L}_i$ denote the matrix representation of the $i$th convolutional layer. Then the equivalent matrix of the entire module is given by the product
\begin{equation}
	\bm{L}=\bm{L}_m\bm{L}_{m-1}\cdots\bm{L}_1.
\end{equation}
For simplicity, we assume circular padding. Together with the kernel symmetry constraints imposed above, each matrix $\bm{L}_i$, $i=1,2,\ldots,m$, is a real symmetric circulant matrix. In general, the product of two symmetric matrices is not necessarily symmetric. However, circulant matrices possess additional algebraic structures that make the analysis more tractable. We recall the following two basic properties.

\begin{lemma}[Commutativity of Circulant Matrices \cite{gray2006toeplitz}]\label{lem:commutative}
	Let $\mathcal{C}_n$ denote the set of all $n\times n$ circulant matrices. For any $\bm{A},\bm{B}\in\mathcal{C}_n$, we have
	\begin{equation}
		\bm{A}\bm{B}=\bm{B}\bm{A}.
	\end{equation}
\end{lemma}

\begin{lemma}[Closure of Circulant Matrices under Multiplication \cite{gray2006toeplitz}]\label{lem:multiply}
	For any $\bm{A},\bm{B}\in\mathcal{C}_n$, their product $\bm{A}\bm{B}$ is also a circulant matrix, namely,
	\begin{equation}
		\bm{A}\bm{B}\in\mathcal{C}_n.
	\end{equation}
\end{lemma}

By Lemma \ref{lem:commutative}, the product of two real symmetric circulant matrices is still symmetric. Indeed, for $\bm{L}_1,\bm{L}_2\in\mathcal{C}_n$ with $\bm{L}_1^T=\bm{L}_1$ and $\bm{L}_2^T=\bm{L}_2$, we have
\begin{equation}
    (\bm{L}_2\bm{L}_1)^T=\bm{L}_1^T\bm{L}_2^T=\bm{L}_1\bm{L}_2=\bm{L}_2\bm{L}_1.
\end{equation}
By induction, this argument extends directly to the product of finitely many real symmetric circulant matrices. Therefore,
\begin{equation}
	\bm{L}=\bm{L}_m\bm{L}_{m-1}\cdots\bm{L}_1
\end{equation}
is symmetric. Moreover, by Lemma \ref{lem:multiply}, the matrix $\bm{L}$ remains circulant. Hence, the equivalent matrix of the multi-layer single-channel convolutional module is again a real symmetric circulant matrix.

To further analyze the definiteness of $\bm{L}$, it remains to understand how the spectrum of a product of circulant matrices is related to the spectra of its individual factors. This relation follows from the spectral representation of circulant matrices, which will be recalled next.

\begin{lemma}[Spectral Decomposition of Real Symmetric Circulant Matrices \cite{gray2006toeplitz}]\label{lem:same_diag}
	Let $\mathcal{SC}_n\subset\mathbb{R}^{n\times n}$ denote the set of all $n\times n$ real symmetric circulant matrices. There exists a real orthogonal matrix $\bm{Q}
    $
    such that every $\bm{A}\in\mathcal{SC}_n$ can be diagonalized as
	\begin{equation}
		\bm{A}=\bm{Q}\bm{\Lambda}_{\bm{A}}\bm{Q}^T,
	\end{equation}
	where $\bm{\Lambda}_{\bm{A}}$ is a real diagonal matrix whose diagonal entries are the eigenvalues of $\bm{A}$.
\end{lemma}

\begin{lemma}[Spectral Property of Products of Simultaneously Diagonalizable Symmetric Matrices]\label{lem:eig_AB}
	Let $\bm{A}$ and $\bm{B}$ be real symmetric matrices. Suppose that there exists a real orthogonal matrix $\bm{Q}
    $ such that $\bm{A}$ and $\bm{B}$ can be simultaneously diagonalized, namely,
	\begin{equation}
		\bm{A}=\bm{Q}\bm{\Lambda}_{\bm{A}}\bm{Q}^T,\qquad \bm{B}=\bm{Q}\bm{\Lambda}_{\bm{B}}\bm{Q}^T,
	\end{equation}
	where $\bm{\Lambda}_{\bm{A}}$, $\bm{\Lambda}_{\bm{B}}$ are the diagonal matrices consisting of the eigenvalues of $\bm{A}$ and $\bm{B}$, respectively. Then their product satisfies
	\begin{equation}
		\bm{A}\bm{B}=\bm{Q}(\bm{\Lambda}_{\bm{A}}\bm{\Lambda}_{\bm{B}})\bm{Q}^T.
	\end{equation}
	Consequently, the eigenvalues of $\bm{A}\bm{B}$ are given by the pointwise products
	\begin{equation}
		\lambda_i(\bm{A}\bm{B})=\lambda_i(\bm{A})\lambda_i(\bm{B}),\qquad i=1,2,\ldots,n,
	\end{equation}
	where $\lambda_i(\bm{A})$ and $\lambda_i(\bm{B})$ are the eigenvalues of $\bm{A}$ and $\bm{B}$ corresponding to the same common eigenvector $\bm{q}_i$, the $i$th column of $\bm{Q}$.
\end{lemma}

\begin{lemma}[Eigenvalue Criterion for Negative Semidefiniteness \cite{horn2012matrix}]\label{lem:negative}
	Let $\bm{A}\in\mathbb{R}^{n\times n}$ be a real symmetric matrix. Then $\bm{A}$ is negative semidefinite, denoted by $\bm{A}\preceq 0$, if and only if all its eigenvalues are non-positive, namely,
	\begin{equation}
		\lambda_i(\bm{A})\leq 0,\qquad i=1,2,\ldots,n.
	\end{equation}
\end{lemma}

By Lemma \ref{lem:negative}, the negative semidefiniteness of $\bm{L}_2\bm{L}_1$ is equivalent to the non-positivity of all its eigenvalues. Since real symmetric circulant matrices share a common real orthogonal eigenbasis by Lemma \ref{lem:same_diag}, Lemma \ref{lem:eig_AB} implies that the eigenvalues of $\bm{L}_2\bm{L}_1$ are given by the pointwise products of the corresponding eigenvalues of $\bm{L}_1$ and $\bm{L}_2$. Therefore, if $\bm{L}_1$ is negative semidefinite and $\bm{L}_2$ is positive semidefinite, then $\bm{L}_2\bm{L}_1$ is negative semidefinite, since
\begin{equation}
	\lambda_i(\bm{L}_2\bm{L}_1)=\lambda_i(\bm{L}_2)\lambda_i(\bm{L}_1)\leq 0,\qquad i=1,2,\ldots,n.
\end{equation}

This argument can be extended to finitely many layers. To construct a negative semidefinite equivalent matrix
\begin{equation}
	\bm{L}=\bm{L}_m\bm{L}_{m-1}\cdots\bm{L}_1,
\end{equation}
a sufficient strategy is to choose $\bm{L}_1$ to be negative semidefinite and require the remaining factors $\bm{L}_2,\bm{L}_3,\ldots,\bm{L}_m$ to be positive semidefinite. Under the simultaneous diagonalization property, the eigenvalues of $\bm{L}$ are given by
\begin{equation}
	\lambda_i(\bm{L})=\prod_{k=1}^{m}\lambda_i(\bm{L}_k),\qquad i=1,2,\ldots,n.
\end{equation}
Hence, if exactly one factor is negative semidefinite and all the remaining factors are positive semidefinite, then all eigenvalues of $\bm{L}$ are non-positive, and consequently $\bm{L}$ is negative semidefinite. In practice, one can follow  the kernel constraints in \eqref{nsd_conditions} and \eqref{psd_conditions} to construct a multi-layer single-channel convolutional module.

When applying PI-DOSnet to the AC equation, 
each convolutional module is required to map a single-channel input to a single-channel output, so that the resulting field can be consistently propagated to the next computational stage. Nevertheless, hidden channels can be introduced within the module to enhance the expressive capacity of the learned spatial operator. This leads to a multi-layer convolutional architecture with one input channel, one output channel, and multiple intermediate channels. Specifically, let the module consist of $m$ convolutional layers, and let $d_\ell$ enote the number of channels in the $\ell$th layer, where
\begin{equation}\label{eq: channel_number_condition}
	d_0=1,\qquad d_m=1,
\end{equation}
while the intermediate channel dimensions $d_1,d_2,\ldots,d_{m-1}$ may vary. Thus, the channel configuration can be represented as $1-d_1-d_2-\cdots-d_{m-1}-1$.

For the $\ell$th convolutional layer, let $\bm{L}^{(\ell)}_{pq}\in\mathbb{R}^{n\times n}$ denote the matrix associated with the convolution kernel from the $q$th input channel of this layer to the $p$th output channel, where
\begin{equation}
	p=1,2,\ldots,d_\ell,\qquad q=1,2,\ldots,d_{\ell-1}.
\end{equation}
If $\bm{y}^{(\ell)}_p$ denotes the $p$th channel after the $\ell$th layer, then the channel-wise update is given by
\begin{equation}
	\bm{y}^{(\ell)}_p=\sum_{q=1}^{d_{\ell-1}}\bm{L}^{(\ell)}_{pq}\bm{y}^{(\ell-1)}_q,\qquad p=1,2,\ldots,d_\ell.
\end{equation}
Since $d_0=d_m=1$, the entire module induces an equivalent matrix $\bm{L}\in\mathbb{R}^{n\times n}$ from the input vector $\bm{y}^{(0)}_1$ to the output vector $\bm{y}^{(m)}_1$.

To describe the equivalent matrix $\bm{L}$, we expand the module along all possible channels from the input channel to the output channel. For a path indexed by
\begin{equation}
	\gamma=(i_1,i_2,\ldots,i_{m-1}),
\end{equation}
where $i_\ell=1,2,\ldots,d_\ell$, define the corresponding pathwise product
\begin{equation}
	\bm{P}_{\gamma}
	=
	\bm{L}^{(m)}_{1i_{m-1}}
	\bm{L}^{(m-1)}_{i_{m-1}i_{m-2}}
	\cdots
	\bm{L}^{(2)}_{i_2i_1}
	\bm{L}^{(1)}_{i_1 1}.
\end{equation}
Then the equivalent matrix of the whole multi-channel convolutional module is
\begin{equation}\label{eq:multi_channel_equiv}
	\bm{L}=\sum_{i_1=1}^{d_1}\sum_{i_2=1}^{d_2}\cdots\sum_{i_{m-1}=1}^{d_{m-1}}\bm{L}^{(m)}_{1i_{m-1}}\bm{L}^{(m-1)}_{i_{m-1}i_{m-2}}\cdots\bm{L}^{(2)}_{i_2i_1}\bm{L}^{(1)}_{i_1 1}=\sum_{\gamma}\bm{P}_{\gamma}.
\end{equation}
In other words, $\bm{L}$ is the sum of the matrix products associated with all input-to-output channel paths.

We recall the following two standard properties, which will be used to analyze the structure of $\bm{L}$.

\begin{lemma}[Vector Space Property of Circulant Matrices \cite{gray2006toeplitz}]\label{lem:add}
	For any $\bm{A},\bm{B}\in\mathcal{C}_n$ and any scalars $\alpha,\beta\in\mathbb{C}$, the linear combination
	\begin{equation}
		\bm{C}=\alpha\bm{A}+\beta\bm{B}
	\end{equation}
	is still a circulant matrix, namely,
	\begin{equation}
		\bm{C}\in\mathcal{C}_n.
	\end{equation}
\end{lemma}

\begin{lemma}[Spectral Property of Linear Combinations of Simultaneously Diagonalizable Symmetric Matrices]\label{lem:eig_AplusB}
	Let $\bm{A},\bm{B}\in\mathbb{R}^{n\times n}$ be real symmetric matrices. Suppose that there exists a real orthogonal matrix $\bm{Q}
    $ such that
	\begin{equation}
		\bm{A}=\bm{Q}\bm{\Lambda}_{\bm{A}}\bm{Q}^{T},\qquad \bm{B}=\bm{Q}\bm{\Lambda}_{\bm{B}}\bm{Q}^{T}, 
	\end{equation}
	where $\bm{\Lambda}_{\bm{A}}$, $\bm{\Lambda}_{\bm{B}}$ are the diagonal matrices consisting of the eigenvalues of $\bm{A}$ and $\bm{B}$, respectively. Then, for any $\alpha,\beta\in\mathbb{R}$, the linear combination $\alpha\bm{A}+\beta\bm{B}$ satisfies
	\begin{equation}
		\alpha\bm{A}+\beta\bm{B}
		=
		\bm{Q}\left(\alpha\bm{\Lambda}_{\bm{A}}+\beta\bm{\Lambda}_{\bm{B}}\right)\bm{Q}^{T}.
	\end{equation}
	Consequently, its eigenvalues are given by
	\begin{equation}
		\lambda_k(\alpha\bm{A}+\beta\bm{B})
		=
		\alpha\lambda_k(\bm{A})+\beta\lambda_k(\bm{B}),
		\qquad k=1,2,\ldots,n,
	\end{equation}
	where $\lambda_k(\bm{A})$ and $\lambda_k(\bm{B})$ correspond to the same common eigenvector.
\end{lemma}

Assume that all channel-to-channel convolution matrices $\bm{L}^{(\ell)}_{pq}$ are real symmetric circulant matrices. By Lemmas \ref{lem:commutative} and \ref{lem:multiply}, each pathwise product $\bm{P}_{\gamma}$ is again a real symmetric circulant matrix. Since all factors along the path are simultaneously diagonalizable, the eigenvalues of $\bm{P}_{\gamma}$ are given by the pointwise products of the corresponding eigenvalues:
\begin{equation}
	\lambda_k(\bm{P}_{\gamma})=\lambda_k\!\left(\bm{L}^{(m)}_{1i_{m-1}}\right)\lambda_k\!\left(\bm{L}^{(m-1)}_{i_{m-1}i_{m-2}}\right)\cdots\lambda_k\!\left(\bm{L}^{(1)}_{i_1 1}\right).
\end{equation}
Each pathwise product $\bm{P}_{\gamma}$ has been shown to be a real symmetric circulant matrix. By Lemma \ref{lem:add}, the sum of these pathwise products in \eqref{eq:multi_channel_equiv} remains circulant. In addition, real symmetry is preserved under finite summation. Therefore, the equivalent matrix $\bm{L}$ is also a real symmetric circulant matrix. The pathwise products share a common real orthogonal eigenbasis, and Lemma \ref{lem:eig_AplusB} yields the following eigenvalue relation:
\begin{equation}
	\lambda_k(\bm{L})=\sum_{\gamma}\lambda_k(\bm{P}_{\gamma}),\qquad k=1,2,\ldots,n,
\end{equation}
where the summation is taken over all input-to-output channel paths $\gamma$.

This observation leads to a simple sufficient condition for the negative semidefiniteness of the whole single-input single-output multi-channel convolutional module. Assume that, for every input-to-output channel path $\gamma$, the corresponding pathwise product $\bm{P}_{\gamma}$ contains an odd number of negative semidefinite factors, while all the remaining factors are positive semidefinite. Then all eigenvalues of $\bm{P}_{\gamma}$ are non-positive, namely,
\begin{equation}
	\lambda_k(\bm{P}_{\gamma})\leq 0,\qquad k=1,2,\ldots,n.
\end{equation}
Thus, each $\bm{P}_{\gamma}$ is negative semidefinite. Since the equivalent matrix $\bm{L}$ is the sum of all pathwise products, the eigenvalue relation above gives
\begin{equation}
	\lambda_k(\bm{L})=\sum_{\gamma}\lambda_k(\bm{P}_{\gamma})\leq 0,\qquad k=1,2,\ldots,n.
\end{equation}
Therefore, $\bm{L}$ is negative semidefinite.

A particularly simple implementation is to impose negative semidefiniteness on all convolution matrices in the first layer and positive semidefiniteness on all convolution matrices in the remaining layers:
\begin{equation}
	\bm{L}^{(1)}_{i_1 1}\preceq 0,\qquad i_1=1,2,\ldots,d_1,
\end{equation}
and
\begin{equation}
	\bm{L}^{(\ell)}_{i_\ell i_{\ell-1}}\succeq 0,\qquad \ell=2,3,\ldots,m.
\end{equation}
Under this design, each channel path contains exactly one negative semidefinite factor, while remaining other factors are positive semidefinite. Consequently, the matrix product associated with each path is negative semidefinite. Since the equivalent matrix $\bm{L}$ is obtained by summing these pathwise products, it follows that \(\bm{L}\) is also negative semidefinite. This result is summarized in the following lemma. 

\begin{lemma}
\label{semidefinite-set-of-neural}
    Let $\mathcal{MCO}$ denote the class of multi-layer convolutional operators $\mathcal{L}_\theta$ satisfying the channel-number constraint \eqref{eq: channel_number_condition}, with the kernel weights in the first layer satisfying \eqref{nsd_conditions} and those in all subsequent layers satisfying \eqref{psd_conditions}. Then for any $\mathcal{L}_{\bm\theta}\in \mathcal{MCO}$, the corresponding equivalent matrix $\bm{L}_{\bm\theta}$ belongs to ${\mathcal{SC}_n}$ and is negative semidefinite.
\end{lemma}

In summary, by explicitly imposing the kernel constraints in \eqref{nsd_conditions} and \eqref{psd_conditions}, we obtain a multi-layer, multi-channel convolutional module whose equivalent operator is negative semidefinite. Consequently, the learned convolutional operator inherits a key structural property required for energy dissipation, providing a foundation for the construction of energy-stable neural operators.

\subsection{Energy Dissipation Preserving PI-DOSnet for linear heat equation}
With a negative semidefinite convolutional module, we now illustrate the construction of an energy-dissipation-preserving PI-DOSnet (EDP-PI-DOSnet) for a linear problem. In the AC equation \eqref{eq:AC}, let us assume $F(u)=0$, resulting in a linear heat equation $\frac{\partial u}{\partial t}=\mathcal{L}u$ with $\mathcal{L}:=\epsilon^2 \Delta$. The energy dissipation for the linear heat equation reads as  
\begin{equation}\label{energy-dissip-linear-eq}
\frac{\mathrm{d}\mathcal{E}(u)}{\mathrm{d}t}=\left\langle u,\frac{\partial u}{\partial t}\right\rangle=\langle u,\mathcal{L}u\rangle\leq 0,
\end{equation}
where 
$
\mathcal{E}(u)=\frac{\epsilon ^2}{2}|\nabla u|^2.
$


The standard PI-DOSnet advances the solution in time through an explicit operator-splitting strategy, which is generally only conditionally stable due to the CFL restriction. More importantly, the resulting scheme does not preserve the energy dissipation law \eqref{energy-dissip-linear-eq}. To address this issue, we modify the linear convolutional component so that the resulting neural operator satisfies the energy dissipation law \eqref{energy-dissip-linear-eq}. Our construction is inspired by the Crank–Nicolson (CN) scheme.
We modify the $i$th block defined in \eqref{PI-DOSnet-ith-block} as
\begin{equation}\label{eq:linear_block}
	\hat{u}^{i}=\psi^{CN}_{dt,\bm{\theta}_i}(\hat{u}^{i-1})=\left(\mathcal{I}-\frac{dt}{2}\mathcal{L}_{\bm{\theta}_i}\right)^{-1}\left(\mathcal{I}+\frac{dt}{2}\mathcal{L}_{\bm{\theta}_i}\right)\hat{u}^{i-1},
\end{equation}
where $\mathcal{L}_{\bm{\theta}_i}$ is a learned approximation of the linear operator $\mathcal{L}$, constructed from convolutional modules whose equivalent matrix satisfies the negative-semidefinite constraint. Since we consider the purely linear case $F(u)=0$, the nonlinear operator is absent. The inverse of $\mathcal{I}-\frac{dt}{2}\mathcal{L}_{\bm{\theta}_i}$ is implemented by a partial pivoting LU decomposition. The update formula \eqref{eq:linear_block} can be viewed as a neural analog of the CN discretization as follows
\begin{equation}
 \frac{\hat{u}^{i+1}-\hat{u}^i}{dt}=\frac 1 2\mathcal{L}_{\bm\theta_i} \left(\hat{u}^{i+1}+\hat{u}^{i}\right).
 \end{equation}
By composing $K$ blocks together, the modified  global mapping in the interval $[0,t]$ is given by
\begin{equation}\label{eq:linear_map}
\hat{u}(t,\bm{x})=\Phi^{CN}_{\bm{\theta}}(u_0,t)(u_0)=\underbrace{\psi^{CN}_{dt,\bm{\theta}_K}\circ\psi^{CN}_{dt,\bm{\theta}_{K-1}}\circ\cdots\circ\psi^{CN}_{dt,\bm{\theta}_1}}_{K\ \mathrm{times}}(u_0),
\end{equation}
where $\bm{x}\in\{\bm{x}_i\}_{i=1}^{N_x}$ denotes a point on the spatial grid. The network output $\hat{u}(t,\bm{x}_j)$ is expected to approximate the exact solution $u(t,\bm{x}_j)$ after training.

Let $V_h\cong\mathbb{R}^n$ be a discrete function space equipped with the discrete $L^2$ inner product $\langle u,v\rangle_h=h\sum_{j=1}^{n}u_jv_j$, where $v_j\approx v(\bm{x}_j)$.
For the linear heat equation, let us introduce the discrete energy as follows
\begin{equation}
	\mathcal{E}_h(u)=\frac{1}{2}\langle \bm{G}_h u,u\rangle_h,
\end{equation}
where $\bm{G}_h\in\mathcal{SC}_n\subset\mathbb{R}^{n\times n}$ is a real symmetric positive semidefinite matrix approximating $-\Delta$, i.e., the traditional finite difference scheme.
We now show that the newly proposed modified PI-DOSnet satisfies the following discrete energy dissipation property.

\begin{proposition}[Discrete Energy Dissipation for the linear heat equation]\label{prop:ed_heat_cn}
 For any $t>0$, $dt$, and ${\mathcal L}_{{\bm\theta}_i}\in\mathcal{MCO}$ with $i=1,\,2,\,\ldots,\,K$, the numerical solution generated by the energy dissipation preserving PI-DOSnet satisfies the discrete energy dissipation property 
 \begin{equation}
	\mathcal{E}_h(\hat{u}^{i})\leq \mathcal{E}_h(\hat{u}^{i-1}), \quad i=1,\,2,\,\ldots,\,K.
\end{equation}



\end{proposition}

\begin{proof}
Let $\bar{u}^{i-1/2}=(\hat{u}^{i}+\hat{u}^{i-1})/2$. From \eqref{eq:linear_block}, we have
\begin{equation}
	\hat{u}^{i}-\hat{u}^{i-1}=dt \mathcal{L}_{\bm{\theta}_i}\bar{u}^{i-1/2}=dt\,\bm{L}_{\bm{\theta}_i}\bar{u}^{i-1/2}.
\end{equation}
Using the symmetry of $\bm{G}_h$, the discrete energy difference can be written as
\begin{equation}
	\mathcal{E}_h(\hat{u}^{i})-\mathcal{E}_h(\hat{u}^{i-1})=\langle \bm{G}_h\bar{u}^{i-1/2},\hat{u}^{i}-\hat{u}^{i-1}\rangle_h=dt\,\langle \bm{G}_h\bar{u}^{i-1/2},\bm{L}_{\bm{\theta}_i}\bar{u}^{i-1/2}\rangle_h.
\end{equation}
Since ${\mathcal L}_{{\bm\theta}_i}\in\mathcal{MCO}$, $\bm{L}_{\bm{\theta}_i}\in\mathcal{SC}_n$ follows from Lemma~\ref{semidefinite-set-of-neural}.
Since $\bm{G}_h\in\mathcal{SC}_n$ and $\bm{L}_{\bm{\theta}_i}\in\mathcal{SC}_n$, according to Lemma~\ref{lem:same_diag}, there exists a real orthogonal matrix $\bm{Q}$ such that
\begin{equation}
	\bm{G}_h=\bm{Q}\bm{\Lambda}_{G}\bm{Q}^{T},\qquad \bm{L}_{\bm{\theta}_i}=\bm{Q}\bm{\Lambda}_{L}\bm{Q}^{T},
\end{equation}
where the diagonal entries satisfy $\lambda_k^{G}\geq 0$ and $\lambda_k^{L}\leq 0$ for $k=1,2,\ldots,n$. Therefore,
\begin{equation}
	\bm{G}_h\bm{L}_{\bm{\theta}_i}=\bm{Q}(\bm{\Lambda}_{G}\bm{\Lambda}_{L})\bm{Q}^{T}
\end{equation}
is negative semidefinite. Hence,
\begin{equation}
	\mathcal{E}_h(\hat{u}^{i})-\mathcal{E}_h(\hat{u}^{i-1})=dt\,\langle \bar{u}^{i-1/2},\bm{G}_h\bm{L}_{\bm{\theta}_i}\bar{u}^{i-1/2}\rangle_h\leq 0.
\end{equation}
This completes the proof.
\end{proof}

\subsection{Energy Dissipation Preserving PI-DOSnet for Nonlinear Gradient Flows}
The observation in the previous subsection naturally motivates the idea of embedding traditional numerical schemes into PI-DOSnet so that the resulting network can inherit their desirable properties. For nonlinear gradient flow problems, i.e. the AC equation, classical energy-stable numerical methods can generally be categorized as explicit, semi-implicit, and fully implicit schemes. However, fully implicit schemes typically require nonlinear iterative solvers, rendering the mapping from the solution at one time step to that at the next highly nonlinear and difficult to characterize explicitly. Consequently, it is challenging to incorporate such schemes into the PI-DOSnet framework while preserving the underlying energy dissipation property.
In this section, we focus on semi-implicit schemes, such as the SAV and R-DVD-1 methods, and integrate them into PI-DOSnet. This construction yields a class of  EDP-PI-DOSnets that inherit the favorable energy stability properties of the underlying numerical schemes.

Inspired by the SAV/CN scheme, we introduce an auxiliary scalar variable $\hat{r}$ and take it as another output of the modified PI-DOSnet. In particular, 
let $\hat{u}^{i-1}, \hat{r}^{i-1}$ and $\hat{u}^{i}, \hat{r}^i$ denote the input and output solution states of the $i$th time block, respectively.
We then design the $i$th block of the network such that
\begin{equation}\label{eq:nonlinear_block}
	\begin{aligned}
		&\frac{\hat{u}^{i}-\hat{u}^{i-1}}{dt}=\mathcal{L}_{\bm{\theta}_i}\hat{u}^{i-1/2}-\hat{r}^{i-1/2}\bar{b}^{i-1/2}:=\bm{L_{\bm\theta_i}}\hat{u}^{i-1/2}-\hat{r}^{i-1/2}\bar{b}^{i-1/2},\\
		&\hat{r}^{i}-\hat{r}^{i-1}=\frac{1}{2}\left\langle\bar{b}^{i-1/2},\hat{u}^{i}-\hat{u}^{i-1}\right\rangle,
	\end{aligned}
\end{equation}
where $\hat{u}^{i}$ denotes the intermediate output of the $i$-th block of EDP-PI-DOSnet, $\mathcal{L}_{\bm{\theta}_i}\in {\mathcal{MCO}}$ is a multi-layer learnable convolutional operator, and
\begin{equation}
	\hat{u}^{i-1/2}=\frac{\hat{u}^{i-1}+\hat{u}^{i}}{2},\qquad \hat{r}^{i-1/2}=\frac{\hat{r}^{i-1}+\hat{r}^{i}}{2}.
\end{equation}
The nonlinear term is evaluated explicitly as
\begin{equation}
	\bar{b}^{i-1/2}=\frac{F'(\bar{u}^{i-1/2})}{\sqrt{\mathcal{E}_1(\bar{u}^{i-1/2})+C_0}},
\end{equation}
where $\bar{u}^{i-1/2}$ is an explicit second-order approximation of $\hat{u}^{i-1/2}$, and $C_0>0$ is chosen so that the quantity under the square root remains positive. For each time block, we use the same explicit half-step predictor to approximate $\hat{u}^{i-1/2}$. Specifically, given the input state $\hat{u}^{i-1}$ of the $i$th time block, the half-step predictor $\bar{u}^{i-1/2}$ can be constructed using different explicit extrapolation methods. One common choice is the explicit Euler half-step predictor, given by
\begin{equation}\label{eq:explicit_half_step_predictor}
    \bar{u}^{i-1/2}=\hat{u}^{i-1}+\frac{dt}{2}\left(\mathcal{L}_{\bm{\theta}_i}\hat{u}^{i-1}-F'\left(\hat{u}^{i-1}\right)\right),\qquad i=1,2,\ldots,K.
\end{equation}
Alternatively, the classical second-order linear extrapolation can be employed, namely,
\begin{equation}\label{eq:second_order_linear_extrapolation}
    \bar{u}^{i-1/2}=\frac{3}{2}\hat{u}^{i-1}-\frac{1}{2}\hat{u}^{i-2},\qquad i=2,3,\ldots,K.
\end{equation}
For the latter approach, the first time block can be initialized using either the explicit Euler half-step predictor or a ghost point that provides an additional historical state $\hat{u}^{-1}$. Both extrapolation methods preserve the discrete energy dissipation law within the proposed framework, although they may yield different approximation accuracies.

\begin{figure}[!htbp]
	\centering
	\includegraphics[width=0.65\linewidth]{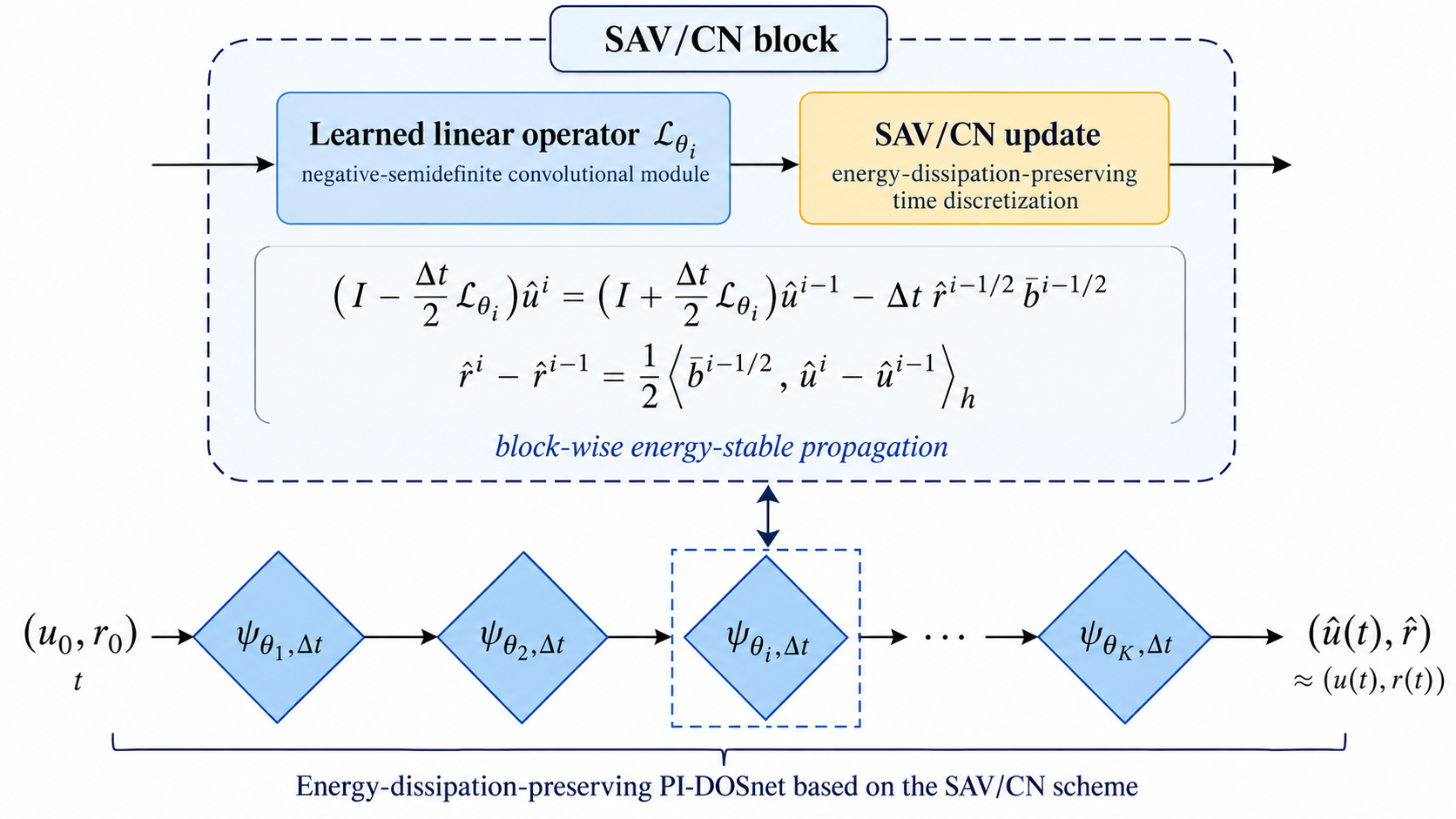}
	\caption{Architecture of the EDP-PI-DOSnet based on the SAV/CN scheme}
	\label{fig:sav_cn_architecture}
\end{figure}

We then rewrite \eqref{eq:nonlinear_block} as
\begin{equation}
\label{eq:nonlinear_block-rewritten}
    \begin{bmatrix}
    \hat{u}^i\\ \hat{r}^i
    \end{bmatrix}=\begin{bmatrix}
    \frac{\mathbf{I}}{dt} -\frac {\bm{L_{\bm\theta_i}}} 2  & * \\ * & 1
    \end{bmatrix}^{-1} \begin{bmatrix}
    \frac{\mathbf{I}}{dt} +\frac {\bm{L}_{\bm{\theta}_i}} 2  & * \\ * & 1
    \end{bmatrix}\begin{bmatrix}
    \hat{u}^{i-1}\\ \hat{r}^{i-1}
    \end{bmatrix},
\end{equation}
which gives an explicit update from the input of the $i$th block to the output of the $i$th block. Therefore, we can construct $i$th block as
\begin{equation}
	\left(\hat{u}^{i},\hat{r}^{i}\right)=\psi^{S-CN}_{dt,\bm{\theta}_i}\left(\hat{u}^{i-1},\hat{r}^{i-1}\right),
\end{equation}
where $\psi_{dt,\bm{\theta}_i}$ corresponds to the product of two matrices in \eqref{eq:nonlinear_block-rewritten}. {We adapt the classical decoupling strategy developed in SAV \cite{shen2018scalar} to solve \eqref{eq:nonlinear_block-rewritten}, with the required matrix inversions performed via LU factorization with partial pivoting.}
For a given time $t$, the modified PI-DOSnet is defined as 
\begin{equation}\label{eq:savcn_map} 
\left(\hat{u}(t,\bm{x}),\hat{r}(t)\right)=\Phi^{S-CN}_{\bm{\theta}}\left(u_0,r_0,t\right)\coloneqq\underbrace{\psi^{S-CN}_{dt,\bm{\theta}_K}\circ\psi^{S-CN}_{dt,\bm{\theta}_{K-1}}\circ\cdots\circ\psi^{S-CN}_{dt,\bm{\theta}_1}}_{K\ \mathrm{times}}\left(u_0,r_0\right). 
\end{equation} 
The network output consists of both the solution component $\hat{u}(t,\bm{x}_j)$ and the auxiliary variable $\hat{r}(t)$, where $\hat{u}(t,\bm{x}_j)$ is expected to approximate the exact solution $u(t,\bm{x}_j)$ after training. The overall architecture of the EDP-PI-DOSnet based on the SAV/CN scheme is illustrated in Fig.~\ref{fig:sav_cn_architecture}.

The above construction is based on the SAV/CN formulation. Alternatively, one can define the nonlinear time-block mapping by using the R-DVD-1 scheme. Since the time-block decomposition and the repeated composition are the same as in the SAV/CN case, we only specify the update rule within each block.
For the \(i\)th time block, the R-DVD-1-based update is given by
\begin{equation}\label{eq:rdvd_block}
	\begin{aligned}
		&\frac{\hat{u}^{i}-\hat{u}^{i-1}}{dt}=\mathcal{L}_{\bm{\theta}_i}\hat{u}^{i-1/2}-F'(\bar{u}^{i-1/2}):=\bm{L}_{\bm{\theta}_i}\hat{u}^{i-1/2}-F'(\bar{u}^{i-1/2}),\\
		&\hat{r}^{i}=\hat{r}^{i-1}+\left\langle F'(\bar{u}^{i-1/2}),\hat{u}^{i}-\hat{u}^{i-1}\right\rangle.
	\end{aligned}
\end{equation}
Here, \(\mathcal{L}_{\bm{\theta}_i}\) is the learned linear operator represented by convolutional modules, and $\hat{u}^{i-1/2}, \bar{u}^{i-1/2}$ are the same in the SAV/CN-based construction.

Compared with the SAV/CN block, the R-DVD-1 block does not require the scalar term \(\hat{r}^{i-1/2}\bar{b}^{i-1/2}\) in the \(u\)-equation. Instead, the nonlinear force \(F'(\bar{u}^{i-1/2})\) is directly incorporated into the update, while the auxiliary variable \(\hat{r}^{i}\) records the corresponding nonlinear energy variation. Therefore, the R-DVD-1 formulation leads to a fully decoupled update for \(\hat{u}^{i}\), followed by an explicit update of \(\hat{r}^{i}\).

We denote the R-DVD-1-based block mapping by
\begin{equation}
	\left(\hat{u}^{i},\hat{r}^{i}\right)=\varphi^{RD-1}_{dt,\bm{\theta}_i}\left(\hat{u}^{i-1},\hat{r}^{i-1}\right).
\end{equation}
The entire mapping is written as
\begin{equation}\label{eq:rdvd_map}
	\left(\hat{u}(t,\bm{x}),\hat{r}(t)\right)=\Phi_{\bm{\theta}}^{RD-1}\left(u_0,r_0,t\right)\coloneqq\underbrace{\varphi^{RD-1}_{dt,\bm{\theta}_K}\circ\varphi^{RD-1}_{dt,\bm{\theta}_{K-1}}\circ\cdots\circ\varphi^{RD-1}_{dt,\bm{\theta}_1}}_{K\ \mathrm{times}}\left(u_0,r_0\right).
\end{equation}
Thus, both the SAV/CN-based and R-DVD-1-based constructions follow the same block-wise operator learning framework, while differing only in the structure-preserving update used inside each time block.

In \eqref{eq:savcn_map} and \eqref{eq:rdvd_map}, the kernel weights of $\mathcal{L}_{{\bm\theta}_i}$ may, in principle, vary across different blocks. However, since $\mathcal{L}_{{\bm\theta}_i}$ is designed to learn the Laplacian operator $\Delta$, which is independent of the time variable $t$, it is natural to share the kernel weights among all blocks. Consequently, all mappings $\mathcal{L}_{{\bm\theta}i}$ are taken to be identical and are denoted uniformly by $\mathcal{L}_{\bm\theta}$. This weight-sharing strategy offers a significant computational advantage: the inverse of $\bm{L}_{\bm\theta}$ needs to be computed only once and can then be reused throughout all blocks. As a result, the computational cost is reduced substantially, leading to improved training efficiency of the proposed network.

For the AC equation, we introduce the discrete modified energy 
\begin{equation}\label{discrete modified energy}
	\widetilde{\mathcal{E}}_{h,m}(u,r)=\frac{1}{2}\langle \bm{G}_h u,u\rangle_h+r^m,
\end{equation}
where $\bm{G}_h\in\mathcal{SC}_n\subset\mathbb{R}^{n\times n}$ is a real symmetric positive semidefinite matrix approximating $-\Delta$, i.e., the traditional finite difference scheme.  
The discrete modified energies corresponding to $m=2$ and $m=1$ are associated with the SAV/CN and R-DVD-1 schemes, respectively.
It should be noted that $\bm{G}_h$ in \eqref{discrete modified energy} is not necessarily identical to $\bm{L}_{{\bm\theta}_i}$. When $\bm{L}_{\bm\theta_i}$ satisfies certain conditions, which will be specified in the following proposition, the $i$th block of the network based on the SAV/CN scheme can be modified as
\begin{equation}\label{eq:savcn_block-modified}
	\begin{aligned}
		&\frac{\hat{u}^{i}-\hat{u}^{i-1}}{dt}=\bm{L}_{\bm{\theta}_i}\hat{u}^{i-1/2}+\hat{r}^{i-1/2}\bar{b}^{i-1/2},\\
		&\hat{r}^{i}-\hat{r}^{i-1}=\frac{1}{2}\left\langle\bm{G}_h\bm{L_{\bm\theta_i}}^\dagger\bar{b}^{i-1/2},\hat{u}^{i}-\hat{u}^{i-1}\right\rangle,
	\end{aligned}
\end{equation}
where $\bm{L}_{\bm\theta_i}^\dagger$ denotes the Moore--Penrose pseudoinverse of $\bm{L}_{\bm\theta_i}$.
This modification ensures the dissipation of the discrete modified energy \eqref{discrete modified energy}.  We next show that the resulting modified PI-DOSnet satisfies the following discrete energy dissipation property.


\begin{proposition}[Unconditional Modified Energy Stability]\label{prop:sav_network_modified_energy}
For any $t>0$, $dt$, and ${\mathcal L}_{{\bm\theta}_i}\in\mathcal{MCO}$ and $\ker{(\bm{G}_h)}=\ker{(\bm{L_{\bm\theta_i}})}$, the numerical solution generated by \eqref{eq:rdvd_block} or \eqref{eq:savcn_block-modified} satisfies the discrete energy dissipation property 
 \begin{equation}
	\widetilde{\mathcal{E}}_{h,m}(\hat{u}^i,\hat{r}^i)- \widetilde{\mathcal{E}}_{h,m}(\hat{u}^{i-1},\hat{r}^{i-1})\leq 0, \quad i=1,\,2,\,\ldots,\,K.
\end{equation}

\end{proposition}

\begin{proof}
We present the proof for $m=2$. The proof for $m=1$ can be completed similarly. 
Using the symmetry of $\bm{G}_{h}$, we have
\begin{equation}
	\frac{1}{2}\langle \bm{G}_{h}\hat{u}^{i},\hat{u}^{i}\rangle_h-\frac{1}{2}\langle \bm{G}_{h}\hat{u}^{i-1},\hat{u}^{i-1}\rangle_h=\left\langle \bm{G}_{h}\hat{u}^{i-1/2},\hat{u}^{i}-\hat{u}^{i-1}\right\rangle_h.
\end{equation}
Moreover, the auxiliary-variable equation in \eqref{eq:savcn_block-modified} gives
\begin{equation}
	(\hat{r}^{i})^2-(\hat{r}^{i-1})^2=\hat{r}^{i-1/2}\left\langle \bm{G}_h\bm{L_{\bm\theta_i}}^\dagger\bar{b}^{i-1/2},\hat{u}^{i}-\hat{u}^{i-1}\right\rangle_h.
\end{equation}
Therefore,
\begin{equation}\label{eq: derivation}
    \widetilde{\mathcal{E}}_{h,m}(\hat{u}^i,\hat{r}^i)- \widetilde{\mathcal{E}}_{h,m}(\hat{u}^{i-1},\hat{r}^{i-1})=\left\langle\bm{G}_h\hat{u}^{i-1/2}+\hat{r}^{i-1/2}\bm{G}_h\bm{L_{\bm\theta_i}}^\dagger\bar{b}^{i-1/2},\hat{u}^{i}-\hat{u}^{i-1}\right\rangle_h.
\end{equation}
Since \(\bm{G}_h,\,\bm{L_{\bm\theta_i}}\in\mathcal{SC}_n\subset\mathbb{R}^{n\times n}\), from Lemma~\ref{lem:same_diag} it follows that $\bm{G}_h=\bm{Q}\bm{\Lambda}_{\bm G}\bm{Q}^T$ and $\bm{L_{\bm\theta_i}}=\bm{Q}\bm{\Lambda}_{\bm L}\bm{Q}^T$. The entries of $\bm{\Lambda}_{\bm G}$ and $\bm{\Lambda}_{\bm L}$ are nonnegative and nonpositive, respectively.
Moreover, since
$\ker{(\bm{G}_h)}=\ker{(\bm{L_{\bm\theta_i}})}$, we have $\bm{G}_h=\bm{G}_h\bm{G}_h^\dagger\bm{G}_h=\bm{G}_h\bm{L_{\bm\theta_i}}^\dagger\bm{L_{\bm\theta_i}}$.
Therefore, \eqref{eq: derivation} can be rewritten as
\begin{equation}
    \begin{aligned}
        \widetilde{\mathcal{E}}_{h,m}(\hat{u}^i,\hat{r}^i)- \widetilde{\mathcal{E}}_{h,m}(\hat{u}^{i-1},\hat{r}^{i-1})&=\left\langle \bm{G}_h \hat{u}^{i-1/2}+\hat{r}^{i-1/2}\bm{G}_h\bm{L_{\bm\theta_i}}^\dagger\bar{b}^{i-1/2},\hat{u}^{i}-\hat{u}^{i-1} \right\rangle_h \\
        &=\left\langle \bm{G}_h\bm{L_{\bm\theta_i}}^\dagger\bm{L_{\bm\theta_i}} \hat{u}^{i-1/2}+\hat{r}^{i-1/2}\bm{G}_h\bm{L_{\bm\theta_i}}^\dagger\bar{b}^{i-1/2},\hat{u}^{i}-\hat{u}^{i-1} \right\rangle_h \\
        &=\left\langle \bm{G}_h\bm{L_{\bm\theta_i}}^\dagger \left(\bm{L_{\bm\theta_i}} \hat{u}^{i-1/2}+ \hat{r}^{i-1/2}\bar{b}^{i-1/2}\right),\hat{u}^{i}-\hat{u}^{i-1} \right\rangle_h\\
        &=\left\langle \bm{G}_h\bm{L_{\bm\theta_i}}^\dagger \left(\frac{\hat{u}^{i}-\hat{u}^{i-1}}{dt}\right),\hat{u}^{i}-\hat{u}^{i-1} \right\rangle_h. 
    \end{aligned}
\end{equation}
Since $\bm{G}_h\bm{L_{\bm\theta_i}}^\dagger=\bm{Q}\bm{\Lambda}_{\bm G}\bm{\Lambda}_{\bm L}^\dagger\bm{Q}^T$, we have $\bm{G}_h\bm{L_{\bm\theta_i}}^\dagger$ is negative semidefinite and
\begin{equation}
	\widetilde{\mathcal{E}}_{h,m}(\hat{u}^i,\hat{r}^i)- \widetilde{\mathcal{E}}_{h,m}(\hat{u}^{i-1},\hat{r}^{i-1})\leq 0.
\end{equation}
This completes the proof.
\end{proof}

\begin{remark} Here, $\bm{G}_h$ denotes the matrix arising from a standard finite difference discretization. Under periodic boundary conditions, its null space is given by \(\ker(\bm{G}_h)=\operatorname{span}\{\mathbf{1}\}\).  
    For \( \bm{L_{\bm\theta_i}}=\operatorname{circ}(a_0,a_1,\ldots,a_{n-1})\in \mathcal{MCO} \), where 
    $a_0=-\sum_{j=1}^{n-1}a_j$, imposing the additional condition \(a_1>0\) ensures that \(\ker(\bm{L_{\bm\theta_i}})=\ker(\bm{G}_h)=\operatorname{span}\{\mathbf{1}\}\). Thus, this condition guarantees that the assumption \(\ker(\bm{L}_{\bm{\theta}_i})=\ker(\bm{G}_h)\) in Theorem~\ref{prop:sav_network_modified_energy}.

    
    
\end{remark}

By Proposition~\ref{prop:sav_network_modified_energy}, the numerical solutions generated by the modified PI-DOSnet satisfy the discrete energy dissipation property
\begin{equation}
	\widetilde{\mathcal{E}}_{h,m}(\hat{u}^{i},\hat{r}^{i})\leq \widetilde{\mathcal{E}}_{h,m}(\hat{u}^{i-1},\hat{r}^{i-1}),\qquad i=1,2,\ldots,K.
\end{equation}
This energy dissipation property is guaranteed directly by the block-wise update mechanism, rather than being enforced indirectly through an additional penalty term in the loss function. Consequently, the proposed architecture has an intrinsic structure-preserving property and achieves unconditional energy stability at the network level. This built-in structure effectively suppresses nonphysical oscillations and artificial energy growth in the prediction of nonlinear gradient flows, without the need for any energy-based regularization or penalty terms in the loss function.

\subsection{Training and Inference Procedures of Energy Dissipation Preserving PI-DOSnet}
As discussed in the previous subsection, the proposed EDP-PI-DOSnet is unconditionally stable with respect to the modified energy. We now describe its training and inference procedures, taking the PI-DOSnet based on the SAV/CN scheme as a representative example. According to \eqref{eq:sav}, for a prediction of PI-DOSnet $(\hat{u},\hat{r})$, the corresponding residuals can be defined as
\begin{equation}\label{eq:residual_sav}
    \begin{aligned}
        &f(t,\bm{x})\coloneqq \hat{u}_t-\epsilon^2\Delta\hat{u}+\frac{\hat{r}(t)}{\sqrt{\mathcal{E}_1(\hat{u})+C_0}}F'(\hat{u}),\\
        &g(t)\coloneqq \hat{r}_t-\frac{1}{2\sqrt{\mathcal{E}_1(\hat{u})+C_0}}\int_{\Omega}F'(\hat{u})\hat{u}_t\,\mathrm{d}\bm{x}.
    \end{aligned}
\end{equation}
Alternatively, the residual can also be constructed directly from the original AC equation, $u_t=\epsilon^2\Delta u-F'(u)$, resulting in
\begin{equation}\label{eq:residual_ac}
    h(t,\bm{x})\coloneqq \hat{u}_t-\epsilon^2\Delta\hat{u}+F'(\hat{u}).
\end{equation}
At the continuous level, the SAV system \eqref{eq:sav} is equivalent to the original AC equation. Consequently, either the pair of residuals $(f,g)$ or the single residual $h$ can be used to construct the physics-informed loss. The former explicitly incorporates the evolution equation for the auxiliary variable, whereas the latter directly imposes the original AC equation. These two formulations provide alternative choices of loss function, allowing users to select the one that best suits their preferences.

Let $\{u_0^{(n)}\}_{n=1}^{N_u}$ denote the set of initial conditions. For each initial condition $u_0^{(n)}$, we sample a set of collection points $\left\{\left(t_{r,j}^{(n)},\bm{x}_{r,j}^{(n)}\right)\right\}_{j=1}^{N_r}$ in the space-time domain. AS in the standard PI-DOSnet framework, the spatial points are selected from the prescribed equispaced mesh points $\{\bm{x}_i\}_{i=1}^{N_x}$, while the temporal points are randomly sampled from the interval $(0,T]$. The residual $f$ is evaluated at the sampled space-time points, and the scalar residual $g$ is evaluated at the corresponding sampled time instances.  In \eqref{eq:residual_sav}, the spatial integral is approximated using the same grid. Temporal derivatives are computed via automatic differentiation with respect to \(t\), whereas spatial derivatives are approximated using numerical differentiation on this grid.

For the SAV-based residual formulation \eqref{eq:residual_sav}, the loss function for the EDP-PI-DOSnet is defined as
\begin{equation}\label{eq:loss_sav}
	\mathfrak{L}(\bm{\theta})=\frac{1}{N_uN_r}\sum_{n=1}^{N_u}\sum_{j=1}^{N_r}\left(\left|f(t_{r,j}^{(n)},\bm{x}_{r,j}^{(n)})\right|^2+\left|g(t_{r,j}^{(n)})\right|^2\right).
\end{equation}
In both the SAV/CN and R-DVD-1 schemes, the quantity $\bar{u}^{i-1/2}$ is introduced as an explicit approximation to the midpoint value $u(t_{i-1}+\Delta t/2,\bm x)$. A common choice is the second-order linear extrapolation given in \eqref{eq:second_order_linear_extrapolation}. Although this extrapolation is simple and computationally efficient, its accuracy may depend on the local temporal evolution of the solution, particularly when relatively large time steps are used. Consequently, the modified energy may deviate from the original energy, and the dissipation of the modified energy guaranteed by EDP-PI-DOSnet does not necessarily imply a monotonic decrease in the original energy at every time step. To address this issue, we propose a learnable extrapolation strategy that adaptively determines the extrapolation coefficients from the available temporal information, rather than relying on manually prescribed coefficients.
Specifically, the midpoint approximation is parameterized as
\begin{equation}
\label{learningablemidpoint}
    \bar{u}^{i-1/2}_{\bm\theta}=\alpha u^{i-1}+(1-\alpha)u^{i-2}+\Delta t\,\mathcal{N}_{\bm\theta}(Z^{i-1}),
\end{equation}
where $\mathcal{N}_{\bm\theta}$ denotes a neural network and $\alpha$ is a learnable extrapolation coefficient. The input of $\mathcal{N}_{\bm\theta}$ is defined as
\begin{equation}
    Z^{i-1}=\left[u^{i-1},u^{i-2},\Delta t\right].
\end{equation}
By optimizing $\alpha$ and the neural network parameters simultaneously, the proposed formulation allows the extrapolation strategy to adaptively select the contribution from previous time levels and provide a nonlinear correction according to the current evolution state.

Since the learnable half-step predictor $\bar{u}^{i-1/2}_{\bm{\theta}}$ depends only on previously computed states, the corresponding nonlinear terms can be evaluated explicitly. Incorporating this predictor into the structure-preserving time discretization schemes preserves the linearity of the resulting systems with respect to the unknowns at the current time step, thereby avoiding nonlinear iterations. Thus, the learnable extrapolation strategy offers greater flexibility while retaining the linear solution structure of the original schemes.
To promote consistency between the evolution of the auxiliary variable $\hat{r}$ and the original energy functional, we introduce the following consistency loss into the total training loss:
\begin{equation}
    \mathcal{L}_{r}=\left|\hat{r}-\sqrt{\mathcal{E}_1(\hat{u})+C_0} \right|^2 .
\end{equation}
By penalizing deviations from this relation, the consistency loss encourages the auxiliary variable to remain consistent with the original energy functional, thereby reducing the discrepancy between the energy dissipation behavior of the modified and original energy functionals.

The initial condition and periodic boundary conditions are enforced through the network architecture. Consequently, no additional loss terms are required. In the proposed framework, the energy dissipation property is incorporated directly into the block-wise update structure of PI-DOSnet, while the residual loss \eqref{eq:loss_sav} drives the predicted solution toward satisfying the governing equations. The overall training procedure of the energy dissipation preserving PI-DOSnet is summarized in Algorithm~\ref{alg:edp_training}.

\begin{algorithm}[H]
	\caption{Training stage of EDP-PI-DOSnet.}
	\label{alg:edp_training}
	\begin{algorithmic}[1]
		\REQUIRE training time interval $[0,T]$, equispaced spatial grid $\{\bm{x}_i\}_{i=1}^{N_x}$, initial function set $\mathfrak{I}=\{u_0^{(n)}\}_{n=1}^{N_u}$, residual point collection $\mathcal{D}_r=\left\{\{(t_{r,j}^{(n)},\bm{x}_{r,j}^{(n)})\}_{j=1}^{N_r}\right\}_{n=1}^{N_u}$, number of time blocks $K$, maximum iteration number $N_{\mathrm{iter}}$, batch size $\zeta$, learning rate $\gamma$.
		\STATE Initialize the trainable parameters $\bm{\theta}$ under the prescribed structural constraints.
		\FOR{$it=1,\ldots,N_{\mathrm{iter}}$}
			\STATE Sample a mini-batch of $\zeta$ initial functions and their corresponding residual points from $\mathfrak{I}$ and $\mathcal{D}_r$.
			\STATE For each sampled initial function $u_0^{(n)}$, set $r_0^{(n)}=\sqrt{\mathcal{E}_1(u_0^{(n)})+C_0}$.
			\STATE Compute the network prediction $\left(\hat{u}(t_{r,j}^{(n)},\bm{x}_{r,j}^{(n)}),\hat{r}(t_{r,j}^{(n)})\right)=\Phi^{S-CN}_{\bm{\theta}}\left(u_0^{(n)},r^0,t_{r,j}^{(n)}\right)$ by applying $K$ energy-stable time blocks.
			\STATE Compute the residual loss $\mathfrak{L}(\bm{\theta})$ using \eqref{eq:loss_sav}.
			\STATE Update the trainable parameters by Adam:
			\begin{equation*}
				\bm{\theta}\leftarrow\bm{\theta}-\gamma\,\nabla_{\bm{\theta}}\mathfrak{L}(\bm{\theta}).
			\end{equation*}
		\ENDFOR
		\ENSURE trained energy dissipation preserving PI-DOSnet solution operator $\Phi_{\bm{\theta}}$.
	\end{algorithmic}
\end{algorithm}

After training, EDP-PI-DOSnet provides an accurate approximation of the local solution operator in the time interval $[0,T]$. For a given initial condition $u_0(\bm{x})$, the trained network $\Phi_{\bm{\theta}_T}$ predicts the solution at any time $t\in[0,T]$ by
\begin{equation}
\hat{u}(t,\bm{x})=\Phi^{S-CN}_{\bm{\theta}_T}(u_0,r_0, t).
\end{equation}
More importantly, EDP-PI-DOSnet retains the long-time prediction capability of the original PI-DOSnet. Using the prediction at time $T$ as the initial condition for the next time interval, the trained network can be recursively applied to approximate the solution over $[T,2T]$, and subsequently over successive time intervals. Due to the energy dissipation preserving property, this recursive procedure enables long-time prediction over arbitrarily large time horizons without retraining the network.


The general inference procedure for EDP-PI-DOSnet that preserves energy
dissipation is summarized in Algorithm \ref{alg:edp_inference}.
It should be noted that the prediction time horizon $T$ need not coincide with the time horizon used during training. For instance, when investigating the temporal evolution of the energy and verifying the energy dissipation property, it is often desirable to obtain predictions at finer time levels. In such cases, a smaller prediction time horizon $T$ can be employed recursively to generate a sequence of solutions that satisfies the discrete energy dissipation law.

\begin{algorithm}[H]
    \caption{Inference stage of EDP-PI-DOSnet.}
    \label{alg:edp_inference}
    \begin{algorithmic}[1]
    \REQUIRE trained network $\Phi_{\bm{\theta}_T}$, initial-condition set $\mathfrak{I}_0$, prediction time length $T$, target time $t^\star$.
    
    \STATE Set $Q=\lfloor t^\star/T\rfloor$ and $\tau=t^\star-QT$.
    \STATE Initialize $\hat{U}_0(\bm{x})=u_0(\bm{x})$ and $\hat{R}_0=r_0$ for any $(u_0(\bm{x}),r_0)\in \mathfrak{I}_0$.
    \FOR{$q=0,1,\ldots,Q-1$}
        \STATE Update $\mathfrak{I}_{q+1}:=\{(\hat{U}_{q+1}(\bm{x}),\hat{R}_{q+1})=\Phi^{S-CN}_{\bm{\theta}_T}(\hat{U}_q(\bm{x}),\hat{R}_q,T) ~\big|~(\hat{U}_q(\bm{x}),\hat{R}_q)\in\mathfrak{I}_{q}\}$.
    \ENDFOR    
        \STATE Predict $(\hat{u}(t^\star,(\bm{x})),\hat{r}(t^\star))=\Phi_{\bm{\theta}_T}(\hat{U}_Q(\bm{x}),\hat{R}_Q,\tau)\in \mathfrak{I}_{t^\star}$ for all $\hat{U}_Q(\bm{x})\in\mathfrak{I}_Q$.
    \ENSURE predicted solution set $\mathfrak{I}_{t^\star}$.
    \end{algorithmic}
\end{algorithm}


\section{Numerical Experiments}\label{sec:experiments}

We validate the energy dissipation preserving property of EDP-PI-DOSnet through a series of numerical experiments on gradient flow and two-phase flow problems. Especially for the AC equation, we compare PI-DOSnet, its penalty function variant, and EDP-PI-DOSnet, and further analyze the spectral properties of their associated convolution operators.

\subsection{Allen--Cahn equation}

In this subsection, we consider the following one-dimensional AC equation with periodic boundary conditions:
\begin{equation}\label{eq4:ac}
    \left\{
    \begin{aligned}
        &\frac{\partial u}{\partial t}-\epsilon^2\Delta u+\gamma(u^3-u)=0,
        && t\in(0,T_{\mathrm{end}}],\quad x\in(-1,1),\\
        &u(0,x)=u_0(x),
        && x\in(-1,1).
    \end{aligned}
    \right.
\end{equation}
Here, $\epsilon>0$ is a parameter controlling the interface width, and $\gamma>0$ is the nonlinear reaction coefficient. We first set $\epsilon=0.02$ and $\gamma=1$. The initial condition is generated by a truncated Fourier series, where the Fourier coefficients are independently sampled from the standard Gaussian distribution $\mathcal{N}(0,1)$. Specifically, we set
\begin{equation}\label{eq4:fourier}
	u_0(x) = 0.1\left(a_0 + \sum_{m=1}^7 \left(a_m\cos(m\pi x) + b_m\sin(m\pi x)\right)\right),
\end{equation}
where $a_0$, $a_m$, and $b_m$ are independent random variables drawn from $\mathcal{N}(0,1)$. The factor $0.1$ is used to control the amplitude of the initial data. The reference solutions are computed using the Chebfun package \cite{driscoll2014chebfun}.

For all simulations, the training end time is set to $T=1$, while the inference end time is set to $T_{\mathrm{end}}=10$. We sample $1{,}000$ initial functions according to \eqref{eq4:fourier}. For each initial function, we select $20$ temporal samples within the training interval and evaluate them on $200$ uniformly distributed spatial grid points, yielding $4{,}000$ collocation points. The standard PI-DOSnet, Penalty PI-DOSnet, and EDP-PI-DOSnet all employ $4$ time blocks. In each time block, the number of convolutional channels is set to $1$-$4$-$1$, and the size of the convolution kernel is chosen as $11$. For PI-DOSnet with penalty function, the weight of the energy loss term in \eqref{eq4:energy_penalty_loss} is set to $\lambda_e=10$. Both the original PI-DOSnet and the energy-penalty PI-DOSnet are retrained to improve their long-time inference performance, with $8{,}000$ epochs used in each training stage. Unless otherwise specified, the reported results are obtained from the final retraining stage. In contrast, EDP-PI-DOSnet is trained only once for $10,000$ epochs, without requiring any subsequent retraining.

\begin{figure}[!htbp]
    \centering
    \begin{subfigure}[b]{0.4\textwidth}
    \centering
    \includegraphics[width=0.8\linewidth]{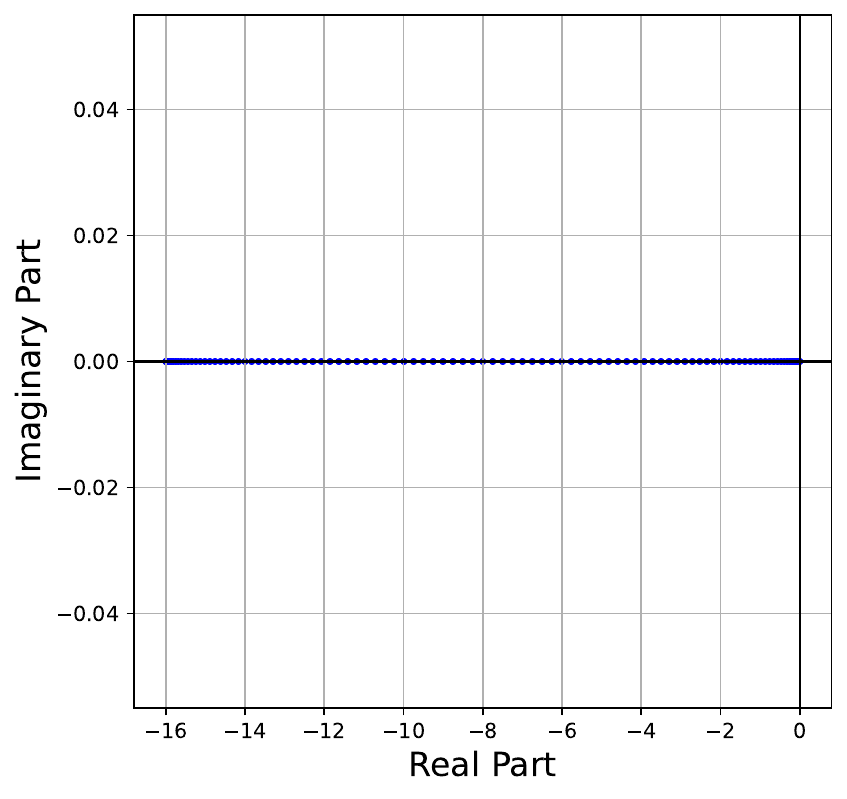}
    \caption{$\Delta_h$}
    \label{fig4:eigvalue-a}
    \end{subfigure}%
    ~
    \begin{subfigure}[b]{0.4\textwidth}
    \centering
    \includegraphics[width=0.78\linewidth]{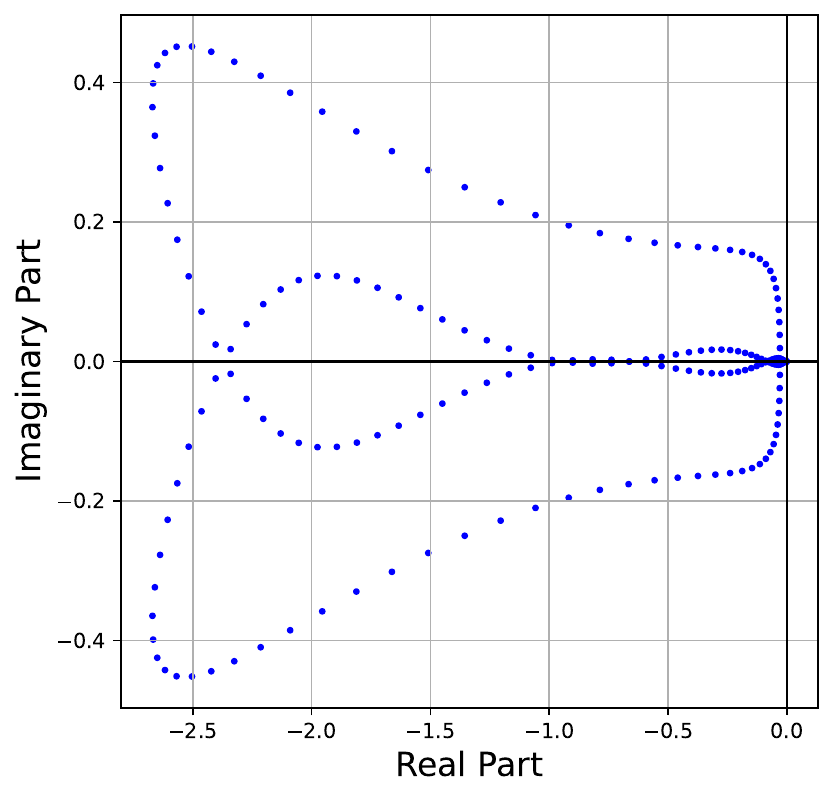}
    \caption{$\bm{L}_{\bm\theta_1}$ in PI-DOSnet}
    \label{fig4:eigvalue-b}
    \end{subfigure}
    \\
    \begin{subfigure}[b]{0.4\textwidth}
    \centering
    \includegraphics[width=0.78\linewidth]{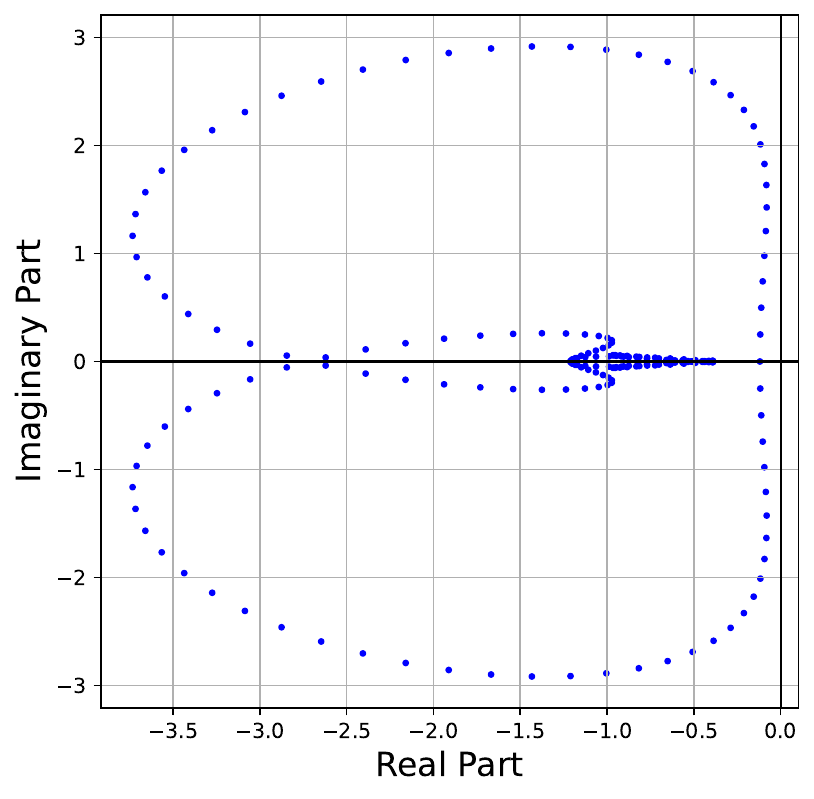}
    \caption{$\bm{L}_{\bm\theta_1}$ in Penalty PI-DOSnet}
    \label{fig4:eigvalue-c}
    \end{subfigure}%
    ~
    \begin{subfigure}[b]{0.4\textwidth}
    \centering
    \includegraphics[width=0.8\linewidth]{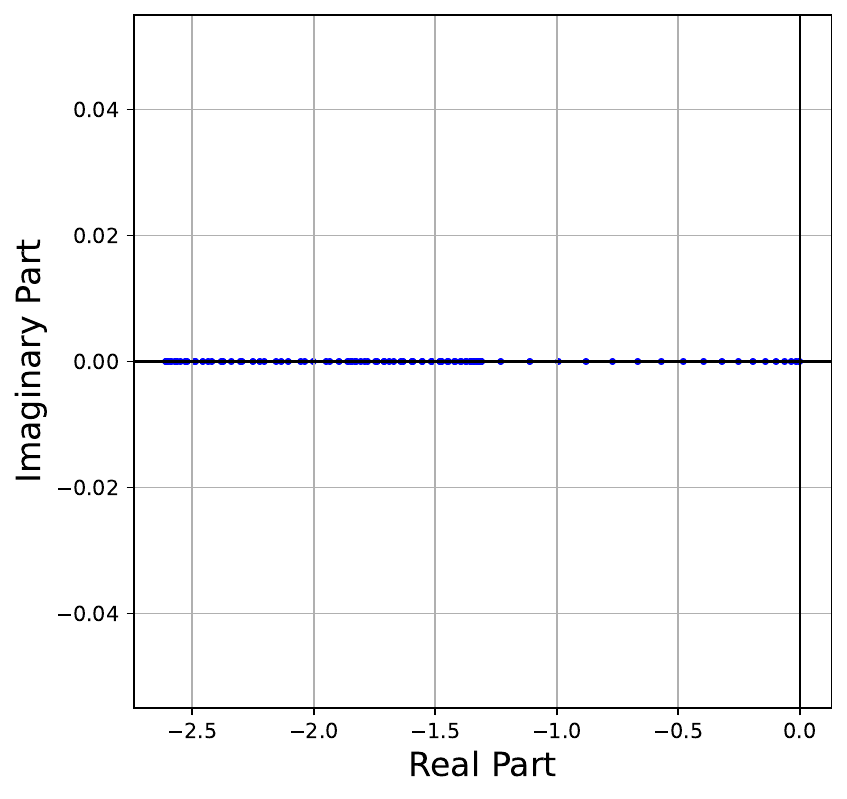}
    \caption{$\bm{L}_{\bm\theta_1}$ in EDP-PI-DOSnet (SAV/CN)}
    \label{fig4:eigvalue-d}
    \end{subfigure}		
    \caption{Eigenvalues distribution of equivalent matrices from different methods. }
    \label{fig4:ac_eigvalue}
\end{figure}

In Fig.~\ref{fig4:ac_eigvalue}, we display the eigenvalue distribution of the learned linear operator ${\cal L}_{\bm\theta}$ for the three PI-DOSnets variants, together with those of the standard second-order central finite difference discretization. For both EDP-PI-DOSnet and the finite difference scheme, all eigenvalues lie on the non-positive real axis, indicating that the corresponding linear operators inherit the negative semidefinite property of $\Delta$. 
In contrast, the eigenvalues of the operators learned by the standard PI-DOSnet and the PI-DOSnet with penalty function are distributed throughout the left half of the complex plane rather than being confined to the non-positive real axis. This suggests that the resulting operators are not negative semidefinite. These results demonstrate that EDP-PI-DOSnet successfully enforces the negative semidefiniteness of the learned linear operator, thereby preserving a key structural property of the underlying differential operator.

\begin{figure}[!htbp]
    \centering
    \includegraphics[width=0.45\textwidth]{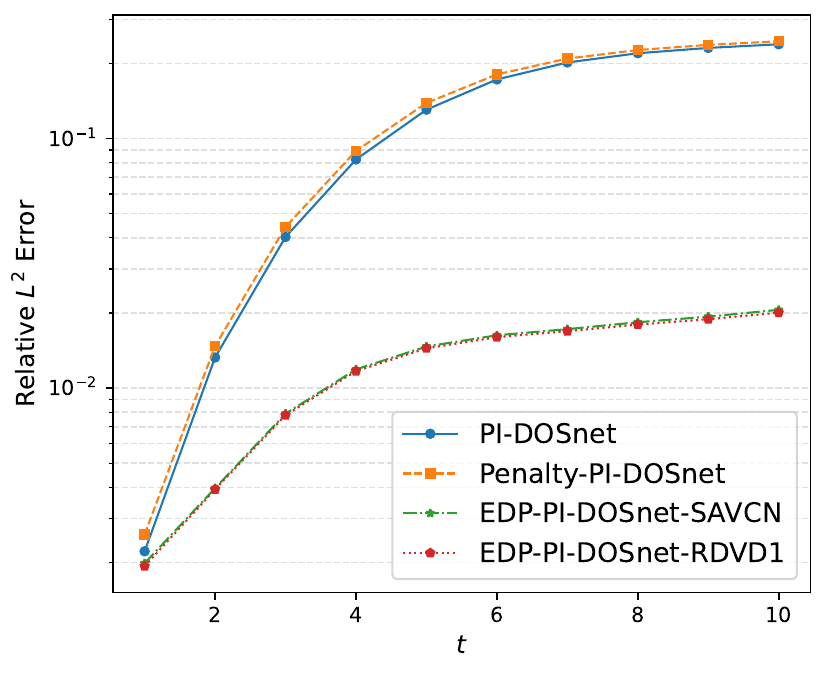}
        \caption{Prediction errors of diffirent PI-DOSnets. 1D AC equation with $T=1$ and $T_{end}=10$.}
    \label{fig4:ac_l2_error}
\end{figure}

After training, all methods are evaluated on a test set of $100$ randomly generated initial functions. Figure~\ref{fig4:ac_l2_error} presents the average prediction error over the test set as a function of time. In this experiment, we adopt a large-step inference strategy with a prediction interval of $T=1$, such that each recursive prediction advances the solution by one time unit. Without retraining, both the standard PI-DOSnet and Penalty PI-DOSnet exhibit relatively large errors in long-time predictions. In contrast, the error curves of the EDP-PI-DOSnet variants based on the R-DVD-1 and SAV/CN schemes nearly overlap, with both achieving substantially smaller errors than the two baseline methods. At $t=10$, the prediction errors of both EDP-PI-DOSnet variants are approximately $2\times10^{-2}$. The performance difference between the two energy-dissipation-preserving methods is negligible. These results demonstrate that EDP-PI-DOSnet achieves accurate and stable long-time prediction without requiring retraining.

\begin{figure}[!htbp]
    \centering
    \begin{subfigure}[b]{0.4\textwidth}
    \centering
    \includegraphics[width=0.8\linewidth]{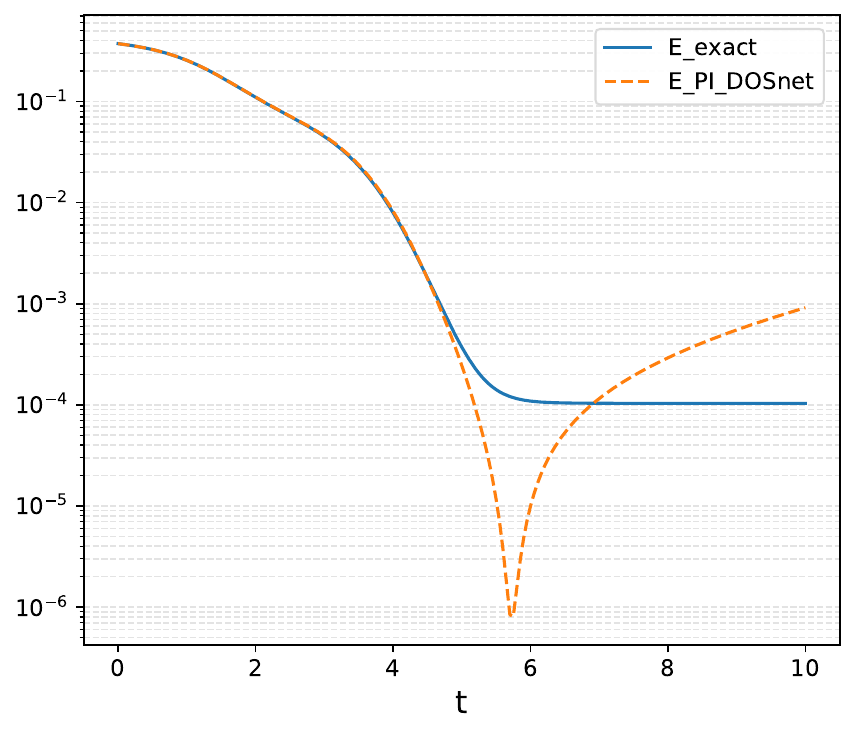}
    \label{fig4:ac_energy_partly-a}
    \end{subfigure}%
    ~
    \begin{subfigure}[b]{0.4\textwidth}
    \centering
    \includegraphics[width=0.8\linewidth]{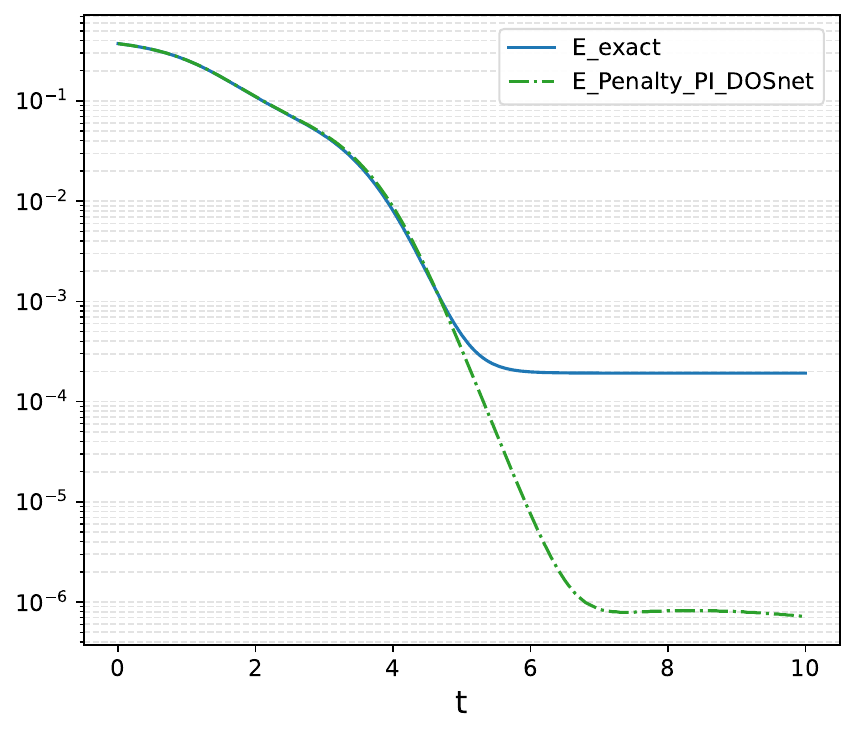}
    \label{fig4:ac_energy_partly-b}
    \end{subfigure}
    \\
    \begin{subfigure}[b]{0.4\textwidth}
    \centering
    \includegraphics[width=0.8\linewidth]{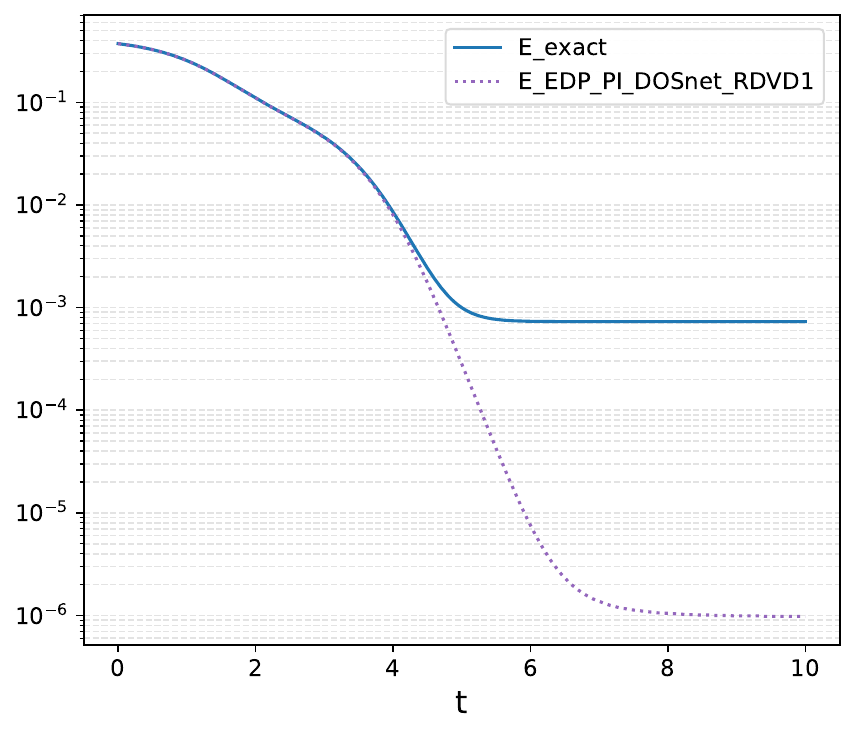}
    \label{fig4:ac_energy_partly-c}
    \end{subfigure}%
    ~
    \begin{subfigure}[b]{0.4\textwidth}
    \centering
    \includegraphics[width=0.8\linewidth]{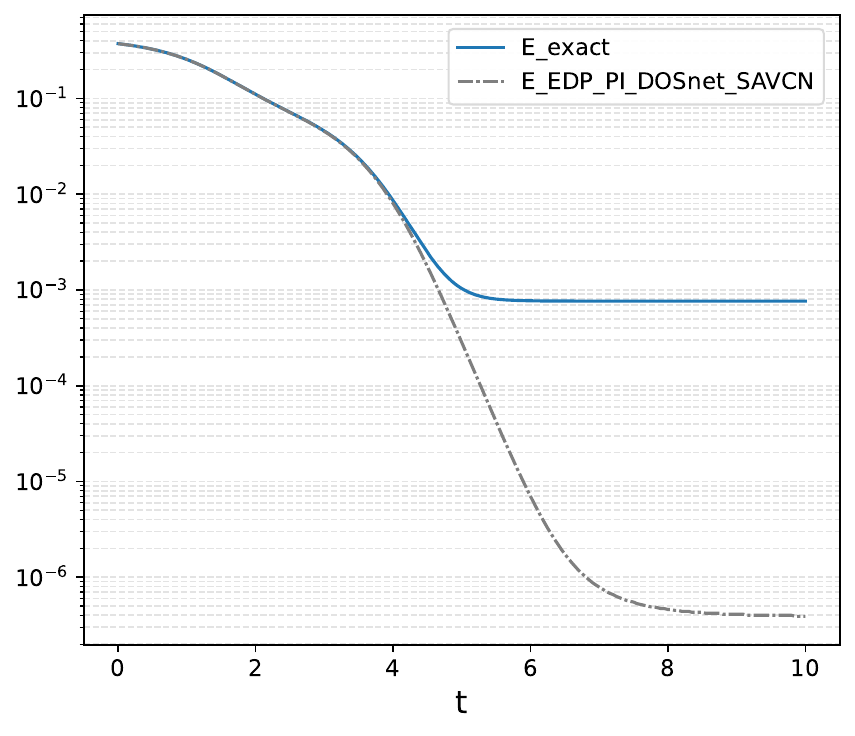}
    \label{fig4:ac_energy_partly-d}
    \end{subfigure}		
    \caption{The modified energy evolution of 1D AC equation (each panel is shifted by a constant determined from its corresponding numerically computed steady-state energy). Here, ``E\_exact'' denotes a curve of numerical modified energies obtained by the Chebfun package with small time step.}
    \label{fig4:ac_energy_partly}
\end{figure}

Next, we study the energy dissipation preserving property of three PI-DOSnets variants. {The discrete energy for the EDP-PI-DOSnet is given by \eqref{discrete modified energy} with $\bm{G}_h={\bm L}_{\bm{\theta}}$.} The discrete energy for the standard PI-DOSnet and the Penalty-PI-DOSnet is computed by \eqref{discrete modified energy} with $\bm{G}_h=\Delta_h$,
where $\Delta_h$ denotes the matrix of the standard second-order central finite difference discretization. 
We set the prediction time length to $\tilde{t}=0.1$ and recursively advance the solution to $T_{\mathrm{end}}$, using the predicted terminal state at each step as the initial condition for the subsequent step. For one random selected initial data, the corresponding energy evolution curves are presented in Fig.~\ref{fig4:ac_energy_partly}. Overall, all methods exhibit a decreasing trend in energy. However, PI-DOSnet exhibits a slight energy increase over $t\in[7,10]$, indicating that energy monotonicity is not strictly preserved. In contrast, both Penalty PI-DOSnet and EDP-PI-DOSnet maintain decreasing (modified) energy profiles. Furthermore, Table~\ref{tab4:count_number_ed} reports the average number of time levels at which the modified energy is non-increasing, evaluated over $100$ test samples for prediction step sizes of $0.05$, $0.1$, and $0.2$. The results show that PI-DOSnet yields the fewest energy non-increasing time levels. Introducing an energy penalty improves this measure by softly encouraging energy dissipation, but does not guarantee strict (modified) energy monotonicity. By contrast, EDP-PI-DOSnet guarantees the non-increase of the modified energy under both the SAV/CN and R-DVD-1 formulations, demonstrating its structure-preserving capability.


\begin{table}[!htbp]
    \caption{Number of time levels satisfying the non-increasing energy condition for the 1D AC equation.}
    \label{tab4:count_number_ed}
    \centering
    \renewcommand{\arraystretch}{1.2}
    \begin{tabular}{cc|cccc} 
        \toprule
        \multicolumn{2}{c|}{Experimental Setting} & \multicolumn{4}{c}{Method} \\ 
        \cmidrule(lr){1-2}\cmidrule(lr){3-6} Time Step & Observed Points & PI-DOSnet & \makecell{Penalty-\\PI-DOSnet} & \makecell{EDP-PI-DOSnet\\(R-DVD-1)} & \makecell{EDP-PI-DOSnet\\(SAV/CN)} \\ \midrule
        0.05  & 200 & 133.50 & 185.45 & 200 & 200\\
        0.1  & 100 & 66.72   & 92.19 & 100 & 100 \\
        0.2  & 50  & 33.39   & 45.47  & 50 & 50 \\
        \bottomrule
    \end{tabular}
\end{table}

Next, we evaluate EDP-PI-DOSnet using pointwise random initial data, where the initial value at each spatial grid point is independently sampled from the uniform distribution $\mathcal{U}(-0.5,0.5)$. Such initial conditions exhibit pronounced high-frequency spatial oscillations and relatively large initial energy, leading to stiff dynamics and posing significant challenges for accurate numerical approximation. This setting therefore provides a demanding test of the robustness and unconditional energy stability of the proposed EDP-PI-DOSnet.
For training over the time interval $(0,1]$, we employ EDP-PI-DOSnet based on the R-DVD-1 scheme with $2$ time blocks. In each block, the convolutional network has a channel configuration of $1$-$4$-$1$, with a kernel size of $5$. For comparison, we also consider a variant of EDP-PI-DOSnet based on the same R-DVD-1 scheme, in which the learnable convolutional operator is replaced by the fixed second-order discrete Laplacian $\Delta_h$. Both models employ the same number of time blocks to ensure a consistent comparison.

\begin{figure}[!htbp]
    \centering
    \begin{subfigure}[b]{0.48\textwidth}
        \centering
        \includegraphics[width=0.8\linewidth]{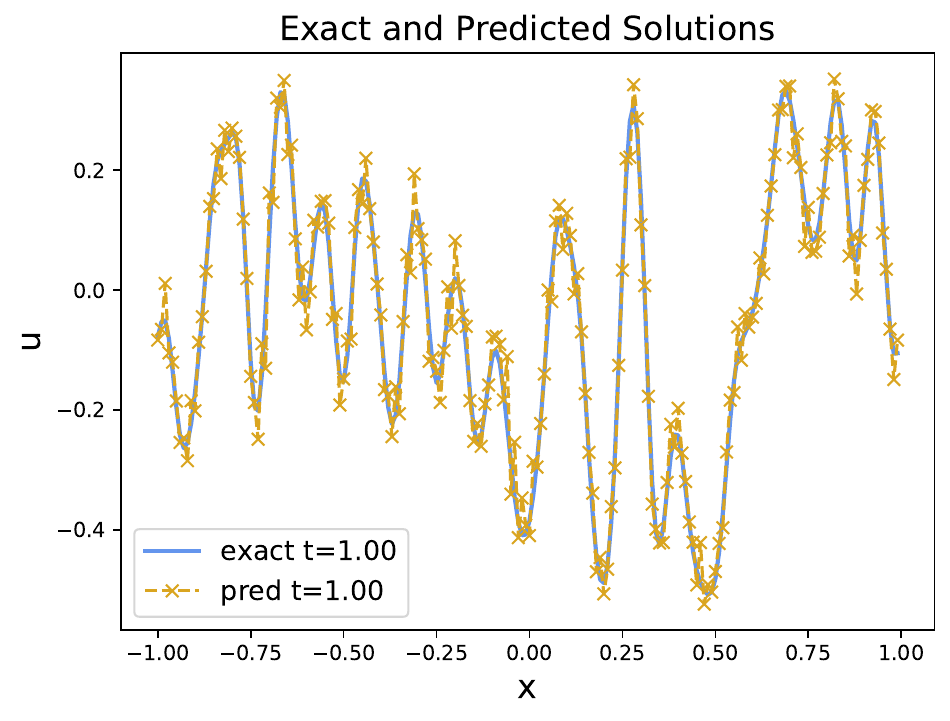}
        \label{fig4:ac_rand_rdvd1-a}
    \end{subfigure}
    \hfill
    \begin{subfigure}[b]{0.48\textwidth}
        \centering
        \includegraphics[width=0.8\linewidth]{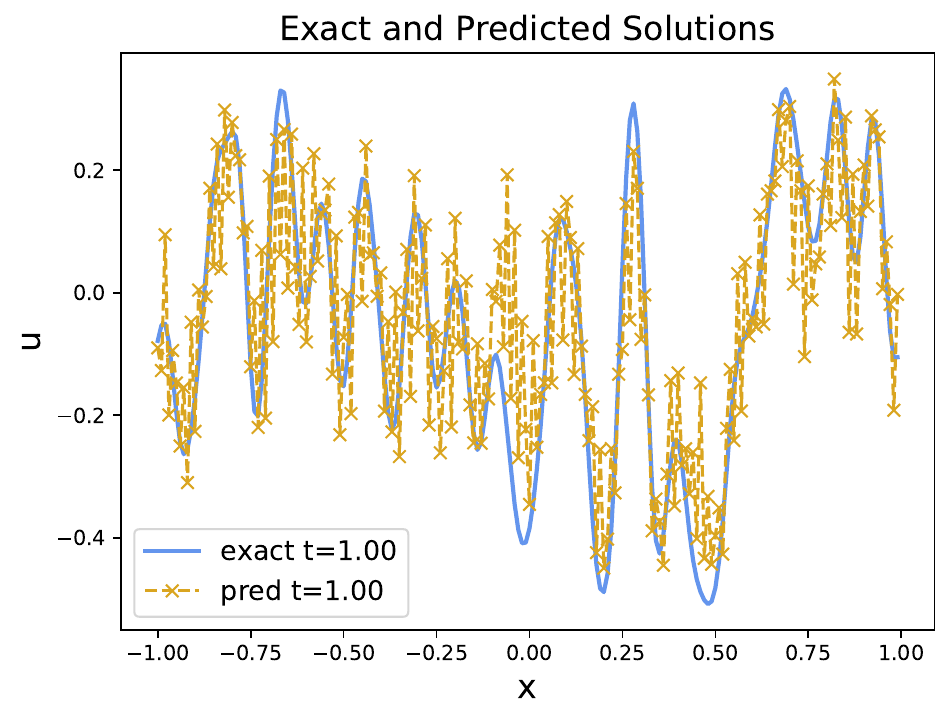}
        \label{fig4:ac_rand_rdvd1-b}
    \end{subfigure}
    \caption{AC equation: reference and predicted solutions at $T=1$ under random initial conditions. Left panel: EDP-PI-DOSnet with a learnable convolutional operator ${\cal L}_{\bm\theta}$; Right panel: EDP-PI-DOSnet with $\Delta_h$.}
    \label{fig4:ac_rand_rdvd1}
\end{figure}

\begin{figure}[!htbp]
    \centering
    \begin{subfigure}[b]{0.48\textwidth}
        \centering
        \includegraphics[width=0.8\linewidth]{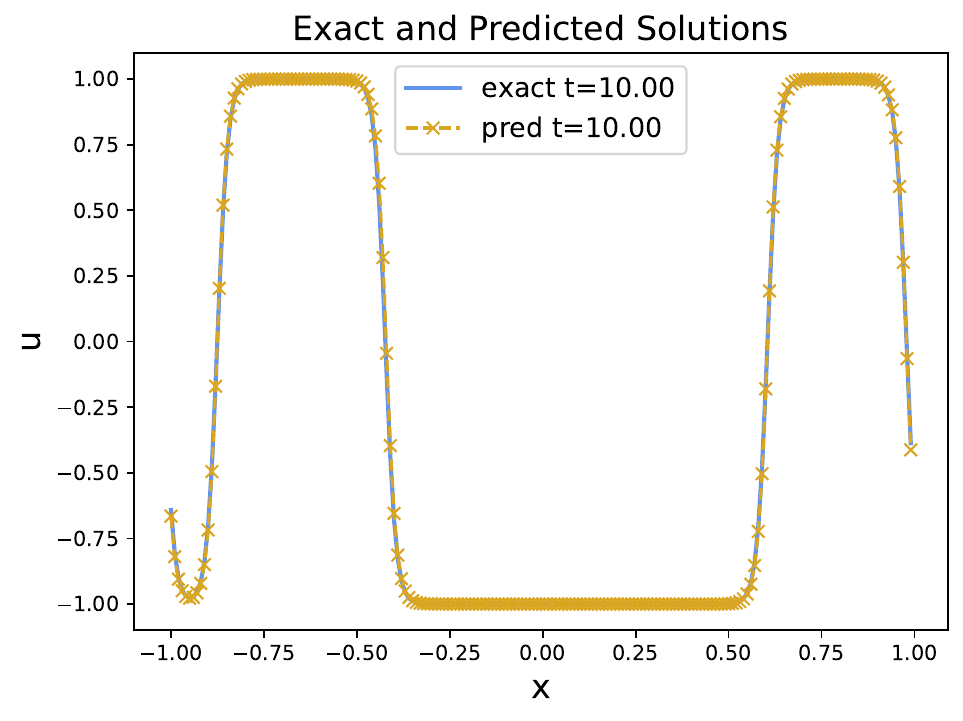}
    \end{subfigure}
    \hfill
    \begin{subfigure}[b]{0.48\textwidth}
        \centering
        \includegraphics[width=0.8\linewidth]{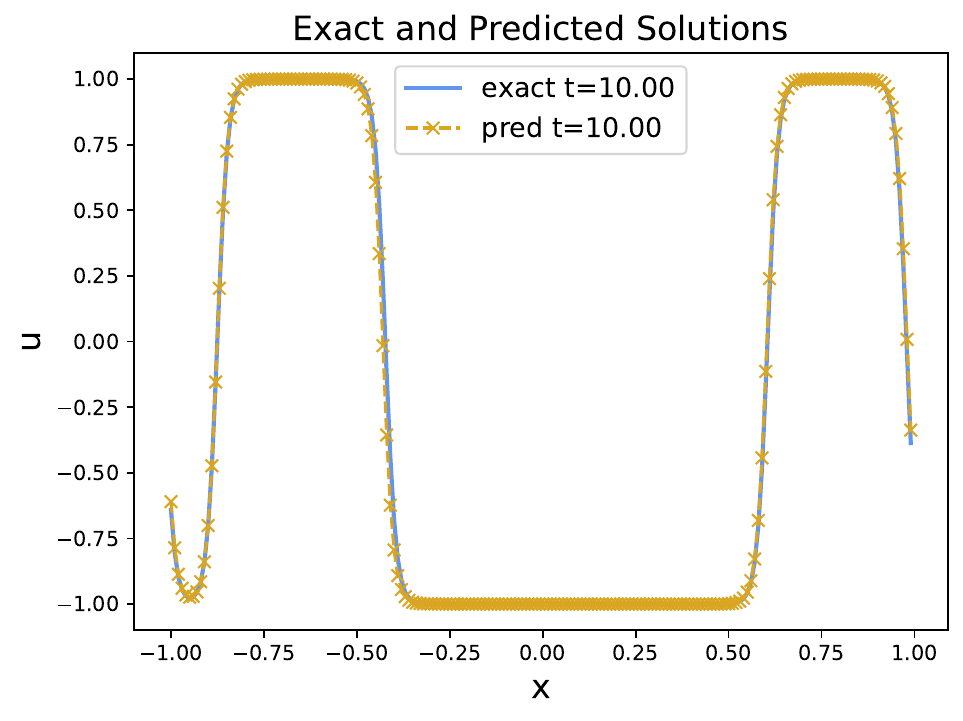}
    \end{subfigure}
    \caption{AC equation: reference and predicted solutions at $T=10$ under random initial conditions. Left panel: EDP-PI-DOSnet with a learnable convolutional operator ${\cal L}_{\bm\theta}$; Right panel: EDP-PI-DOSnet with $\Delta_h$.}
    \label{fig4:ac_rand_rdvd1_Tend}
\end{figure}

Figures~\ref{fig4:ac_rand_rdvd1} and~\ref{fig4:ac_rand_rdvd1_Tend} present the solutions predicted by the two methods at $t=1$ and $t=10$, respectively. At $t=1$, EDP-PI-DOSnet with a learnable convolutional operator ${\cal L}_{\bm\theta}$ shows good agreement with the reference solution, particularly in capturing the main spatial structures arising from the highly oscillatory random initial data. The corresponding prediction errors are $1.531\times10^{-1}$ and $4.854\times10^{-2}$ at $t=1$ and $t=10$, respectively. In comparison, EDP-PI-DOSnet with the fixed discrete Laplacian $\Delta_h$ yields errors of $5.291\times10^{-1}$ and $3.643\times10^{-2}$ at the same time instances. These results demonstrate that the learnable convolutional operator significantly improves prediction accuracy at $t=1$, while the fixed discrete Laplacian achieves a slightly smaller error at $t=10$. Overall, both variants exhibit stable long-time prediction behavior, highlighting the robustness of EDP-PI-DOSnet when applied to highly oscillatory random initial data.

\begin{figure}[!htbp]
    \centering
    \begin{subfigure}[b]{0.48\textwidth}
        \centering
        \includegraphics[width=0.8\linewidth]{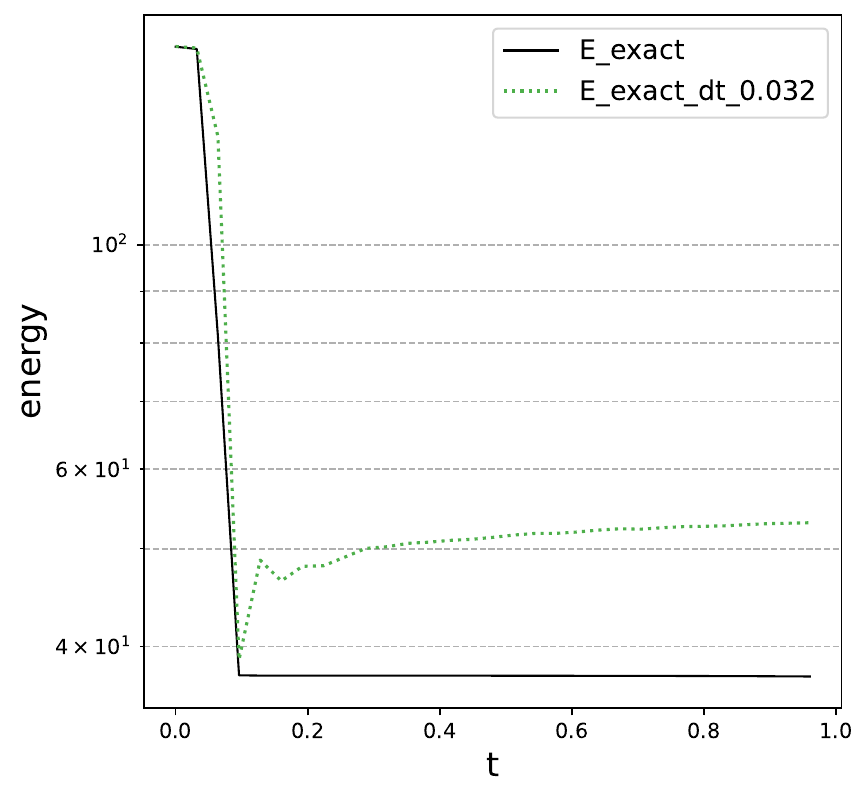}
    \end{subfigure}
    \hfill
    \begin{subfigure}[b]{0.48\textwidth}
        \centering
        \includegraphics[width=0.8\linewidth]{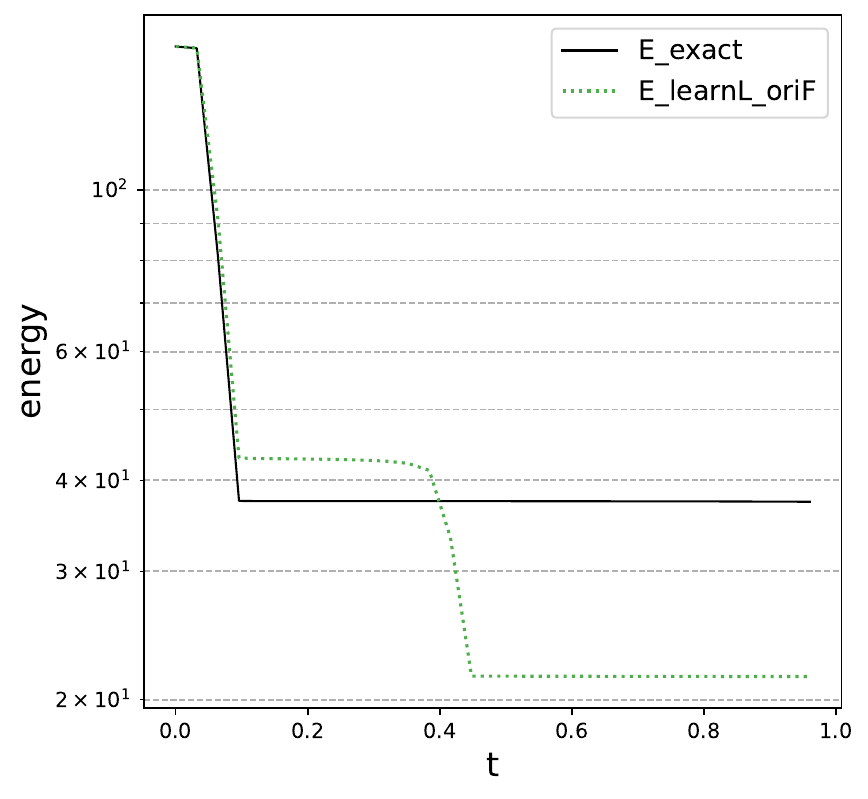}
    \end{subfigure}
    \caption{Comparison of the original energy evolution obtained by different extrapolation strategies.}
    \label{fig: ac_energy_comparison}
\end{figure}

As discussed earlier, EDP-PI-DOSnet guarantees the dissipation of the modified energy, but not necessarily that of the original energy. We therefore conduct a numerical experiment to investigate the dissipation behavior of the original energy and assess whether the learnable midpoint approximation \eqref{learningablemidpoint} can improve its preservation.
This experiment is conducted on a uniform spatial grid with $N_x=100$ points. The initial condition is generated by independently sampling the value at each spatial grid point from the uniform distribution $\mathcal{U}(-0.02,0.02)$. The parameters of the AC equation are set to $\epsilon=1$ and $\gamma=100$. The learnable extrapolation model is trained over the time interval $[0,0.032]$, which is divided into 4 blocks. The trained model is then recursively applied to predict the solution up to $T_{\mathrm{end}}=0.96$, allowing us to assess its long-time behavior. 
Fig.~\ref{fig: ac_energy_comparison} compares the evolution of the original energy obtained using two different midpoint approximations.
The left panel corresponds to the second-order extrapolation $\bar{u}^{i-1/2}=\frac{3}{2}u^{i-1}-\frac{1}{2}u^{i-2}$, while the right panel presents the results obtained using the proposed learnable extrapolation strategy. The second-order extrapolation exhibits noticeable fluctuations in the original energy during long-time evolution, indicating that its monotonic decay is not strictly preserved. In contrast, the learnable extrapolation produces a smoother energy evolution with a more consistent decreasing trend. These results suggest that the proposed learnable extrapolation strategy can improve the dissipation behavior of the original energy during long-time prediction, although strict monotonicity is not theoretically guaranteed.

\subsection{Mass-Conserved Allen--Cahn--Navier--Stokes Equations}
\label{subsec:mass-conserved-ac-ns}

As a more complex application of EDP-PI-DOSnet, we consider the mass-conserved Allen--Cahn--Navier--Stokes (CAC--NS) system for incompressible  two-phase flows. The phase-field variable $\phi$ distinguishes the two fluid phases, with $\phi \approx 1$ and $\phi \approx -1$ corresponding to the two bulk phases, while the transition region represents the diffuse interface between them. The velocity, pressure, and chemical potential are denoted by $\bm{u}$, $p$, and $\mu$, respectively.
Let $\Omega \subset \mathbb{R}^{2}$ be a bounded domain. The mass correction is
introduced through the projected chemical potential $\mu-\overline{\mu}$, where
\begin{equation}
\label{eq:mu-average}
  \overline{\mu}
  = \frac{1}{|\Omega|}\int_{\Omega}\mu\,\mathrm{d}\bm{x}
\end{equation}
is the spatial average of $\mu$. The CAC--NS equations are
\begin{subequations}
\label{eq:acns-continuous}
\begin{gather}
  \phi_t + \nabla \cdot (\bm{u}\phi)
  = -M\left(\mu - \overline{\mu}\right),
  \label{eq:acns-phase} \\
  \mu
  = \lambda\left(-\Delta \phi + f(\phi)\right),
  \label{eq:acns-chemical} \\
  \bm{u}_t + (\bm{u}\cdot\nabla)\bm{u}
  - \nu \Delta \bm{u} + \nabla p + \phi\nabla\mu
  = \bm{0},
  \label{eq:acns-momentum} \\
  \nabla\cdot\bm{u} = 0.
  \label{eq:acns-incompressibility}
\end{gather}
\end{subequations}
Here \(M>0\) is the mobility, \(\nu>0\) is the kinematic viscosity, \(\lambda>0\) is the interfacial energy coefficient, and \(f(\phi)=F'(\phi)\), with \(F(\phi)=\frac{(\phi^2-1)^2}{4\epsilon^2}\) denoting double-well potential. The projected term \(\mu-\overline{\mu}\) enforces phase-mass conservation. Under periodic or compatible no-flux boundary conditions, the following mass conservation law holds:
\begin{equation}
\label{eq:phase-mass-conservation}
    \frac{\mathrm{d}}{\mathrm{d}t}\int_{\Omega}\phi\,\mathrm{d}\bm{x}=0.
\end{equation}

The energy functional for the CAC--NS system is:
\begin{equation}
\label{eq:acns-energy}
  \mathcal{E}(\bm{u},\phi)
  =
  \int_{\Omega}
  \frac{1}{2}|\bm{u}|^2
  + \lambda\left(\frac{1}{2}|\nabla\phi|^2 + F(\phi)\right)
  \,\mathrm{d}\bm{x},
\end{equation}
which satisfies the energy dissipation law
\begin{equation}
\label{eq:acns-energy-law}
  \frac{\mathrm{d}}{\mathrm{d}t}\mathcal{E}(\bm{u},\phi)
  =
  -\nu\|\nabla\bm{u}\|_{L^2(\Omega)}^2
  -M\|\mu-\overline{\mu}\|_{L^2(\Omega)}^2
  \le 0.
\end{equation}

The network architecture used to treat the CAC--NS system is inspired by the DIEQ scheme proposed in \cite{yang2021novel}. DIEQ is closely related in spirit to IEQ and SAV methods and belongs to the class of auxiliary-variable-based energy-stable schemes. In the CAC--NS setting, its key idea is to combine the IEQ treatment of the nonlinear double-well potential with a specially designed nonlocal auxiliary variable for the nonlinear coupling terms. The latter exploits the zero-energy-contribution property of the advection and capillary-stress terms, which enables the coupled system to be reformulated into an equivalent form suitable for explicit treatment of the coupling terms. As a result, the resulting time discretization is linear, fully decoupled, and equipped with a provable discrete energy dissipation law.

We embed one DIEQ-type time step into the network as a structure-preserving block, namely,
\begin{equation}
	(\phi^{n+1},\bm{u}^{n+1},p^{n+1})=\mathcal{N}_{\bm\theta}^{\mathrm{DIEQ}}(\phi^{n},\bm{u}^{n},p^{n}).
\end{equation}
This block follows the algebraic update structure of the DIEQ scheme, while the phase-field Laplace operator is relaxed to a trainable convolutional operator $\mathcal{L}_{\bm\theta}$.
The convolutional operator $\mathcal{L}_{\bm\theta}$ is constrained to be negative semidefinite and to have zero column sum. The conditions and construction ensuring the negative semidefiniteness of this operator have been discussed in detail in the preceding section and are therefore omitted here. The zero-column-sum condition is essential for the mass-conserved phase equation, since it guarantees that the contribution of the relaxed operator vanishes after spatial integration,
\begin{equation}
	(1,\mathcal{L}_{\bm\theta}\phi)=0.
\end{equation}
Consequently, the relaxed operator does not change the spatial average of the phase variable. In practice, this constraint is imposed in the first convolutional layer, chosen as the negative-semidefinite layer, by taking the center kernel coefficient to be negative and all surrounding coefficients to be nonnegative, together with the zero-sum constraint on the kernel.

The proposed network is evaluated on two-droplet coalescence and coarsening problems in the periodic domain $\Omega=[0,L]^2$, with $L=2\pi$. Periodic boundary conditions are imposed in both spatial directions, and the initial velocity and pressure are set to $\bm{u}_0=\bm{0}$ and $p_0=0$, respectively. To construct the initial phase field, we introduce the diffuse droplet profile
\begin{equation}
    \eta_r(\bm{x};\bm{c})
    =\frac{1}{2}\left[1+\tanh\left(\frac{r-d_{\mathrm{per}}(\bm{x},\bm{c})}{1.5\epsilon}\right)\right],
\end{equation}
where $r$ and $\bm{c}$ denote the droplet radius and center, respectively, and $d_{\mathrm{per}}$ is the shortest Euclidean distance under periodic boundary conditions. The initial phase field is defined as a smooth union of two droplets:
\begin{equation}
    \phi_0(\bm{x})
    =s\left[2\left(1-\eta_{r_1}(\bm{x};\bm{c}_1)\right)
    \left(1-\eta_{r_2}(\bm{x};\bm{c}_2)\right)-1\right],
\end{equation}
where the background phase $s$ takes the values $-1$ and $1$ with equal probability. The droplet radii are specified differently for the two configurations. For the coalescence problem, the two droplets have equal radii uniformly sampled from $[0.55,0.70]$ and are initially separated by a small gap. For the coarsening problem, the droplets have unequal radii, with the larger radius sampled uniformly from $[0.19L,0.23L]$ and the smaller radius chosen to be approximately $40\%$--$55\%$ of the larger one. In both cases, the droplet centers and orientations are randomly selected, taking periodicity into account.

The spatial domain is discretized using a uniform $90\times90$ periodic grid. The physical parameters are set to $\epsilon=0.1$, $\lambda=0.01$, $M=50$, $\nu=1$, and $B=1$, with the stabilization parameter $S=2$.
The network consists of $20$ DIEQ blocks and employs a shared trainable convolutional operator composed of two linear convolutional layers with $7\times7$ kernels, four hidden channels, and circular padding. The training dataset contains $400$ initial conditions, including $200$ coalescence and $200$ coarsening realizations. The network is trained for $8,000$ epochs using the Adam optimizer. 
The network is trained by minimizing a physics-informed loss derived from the CAC--NS system \eqref{eq:acns-continuous}. Given the network prediction $(\phi_{\theta},\bm{u}_{\theta},p_{\theta})$, the corresponding chemical potential is computed as
\begin{equation}
	\mu_{\theta}=\lambda\left(-\Delta\phi_{\theta}+f(\phi_{\theta})\right),\qquad \overline{\mu}_{\theta}=\frac{1}{|\Omega|}\int_{\Omega}\mu_{\theta}\,\mathrm{d}\bm{x}.
\end{equation}
The corresponding phase, momentum residuals are
\begin{equation}
	\label{eq:cacns-physical-residuals}
	\begin{aligned}
        &\mathcal{R}_{\phi,\theta}=\partial_t\phi_{\theta}+\nabla\cdot(\bm{u}_{\theta}\phi_{\theta})+M(\mu_{\theta}-\overline{\mu}_{\theta}),\\
		&\mathcal{R}_{\bm{u},\theta}=\partial_t\bm{u}_{\theta}+(\bm{u}_{\theta}\cdot\nabla)\bm{u}_{\theta}-\nu\Delta\bm{u}_{\theta}+\nabla p_{\theta}+\phi_{\theta}\nabla\mu_{\theta}.
	\end{aligned}
\end{equation}
The physical loss is then given by
\begin{equation}
	\label{eq:cacns-physical-loss}
	\mathfrak{L}_{\mathrm{phys}}(\theta)=\|\mathcal{R}_{\phi,\theta}\|_2^2+\omega_{\mathrm{momentum}}\|\mathcal{R}_{\bm{u},\theta}\|_2^2.
\end{equation}
After training, the learned evolution map is first applied over the interval $[0,0.25]$ and subsequently reapplied recursively to extend the prediction to $T_{\mathrm{end}}=0.5$. The network is evaluated on $20$ independent test initial conditions, consisting of $10$ coalescence and $10$ coarsening realizations.

\begin{figure}[!htbp]
    \centering
    \includegraphics[width=0.8\textwidth]{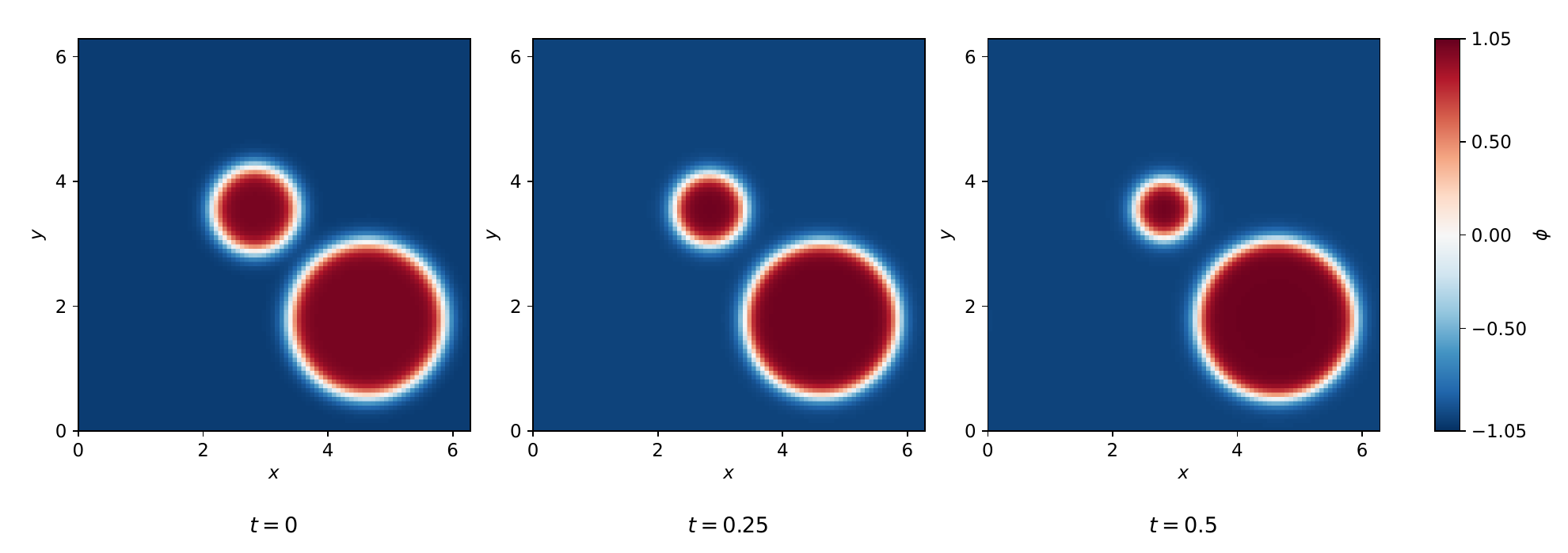}
    \vspace{0.3cm}
    \includegraphics[width=0.8\textwidth]{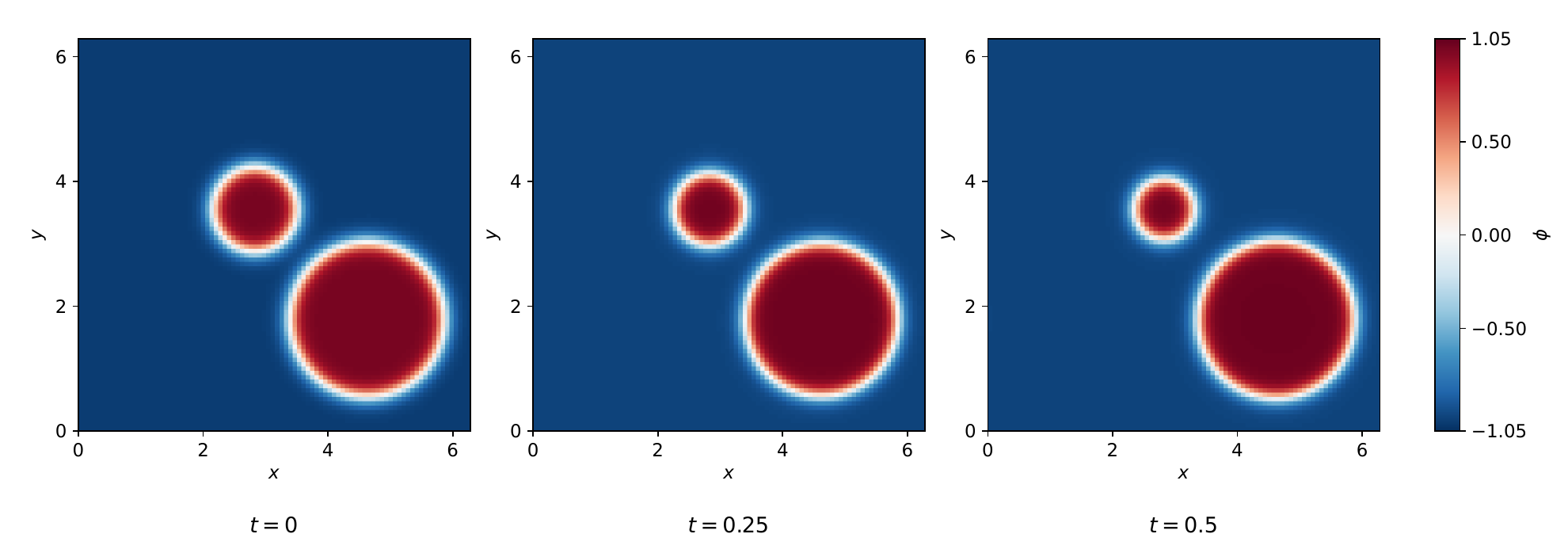}
    \caption{Comparison of the reference and predicted phase-field evolution for the two-droplet coarsening problem. The top and bottom panels show the reference and predicted solutions, respectively.}
    \label{fig:cacns-coarsening}
\end{figure}

\begin{figure}[!htbp]
    \centering
    \includegraphics[width=0.8\textwidth]{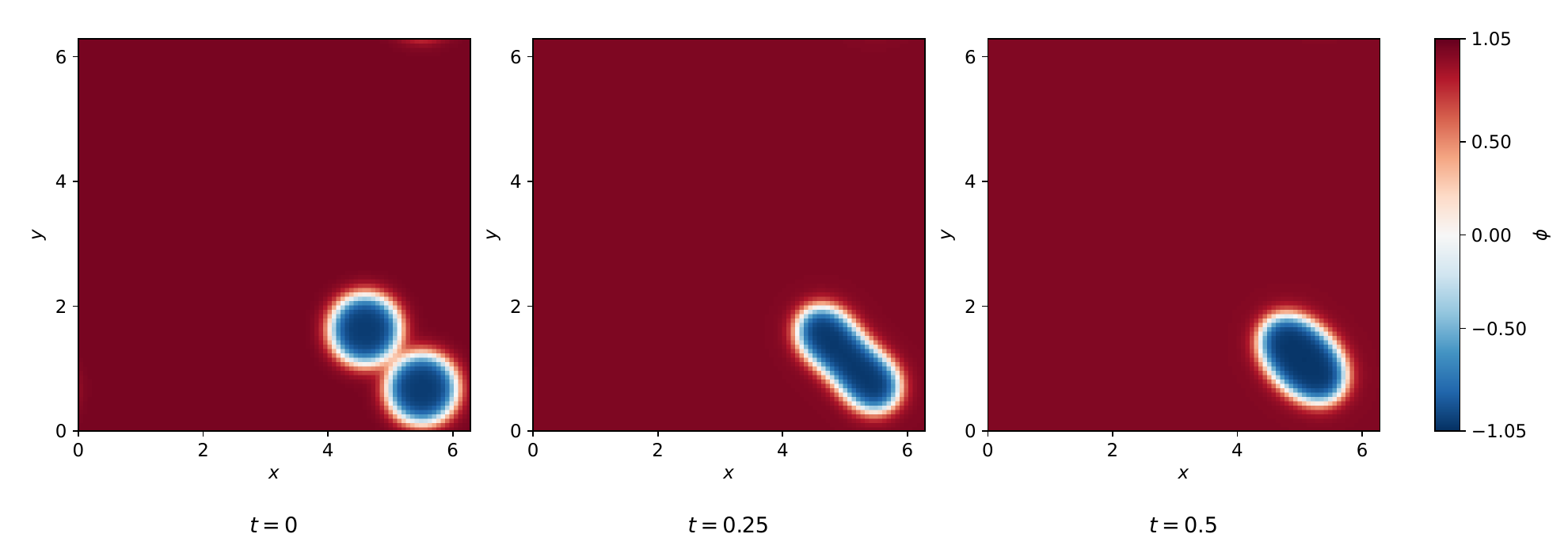}
    \vspace{0.3cm}
    \includegraphics[width=0.8\textwidth]{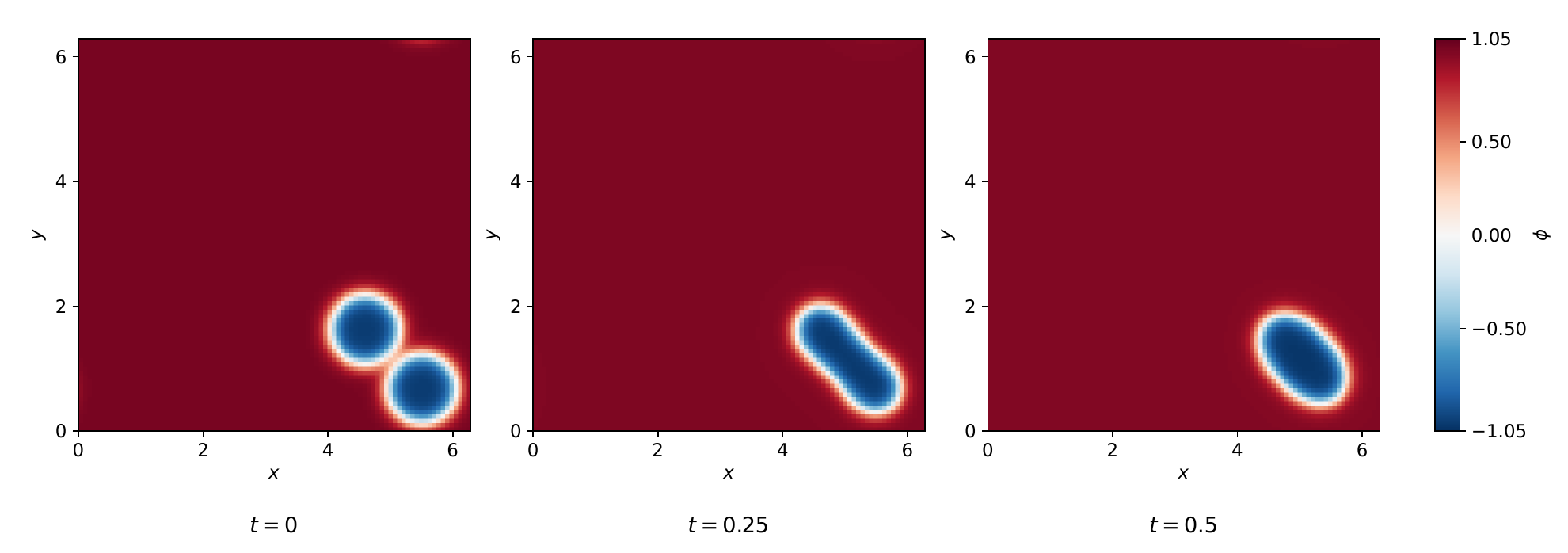}
    \caption{Comparison of the reference and predicted phase-field evolution for the two-droplet coalescence problem. The top and bottom panels show the reference and predicted solutions, respectively.}
    \label{fig:cacns-coalescence}
\end{figure}

\begin{figure}[!htbp]
    \centering
    \includegraphics[width=0.7\textwidth]{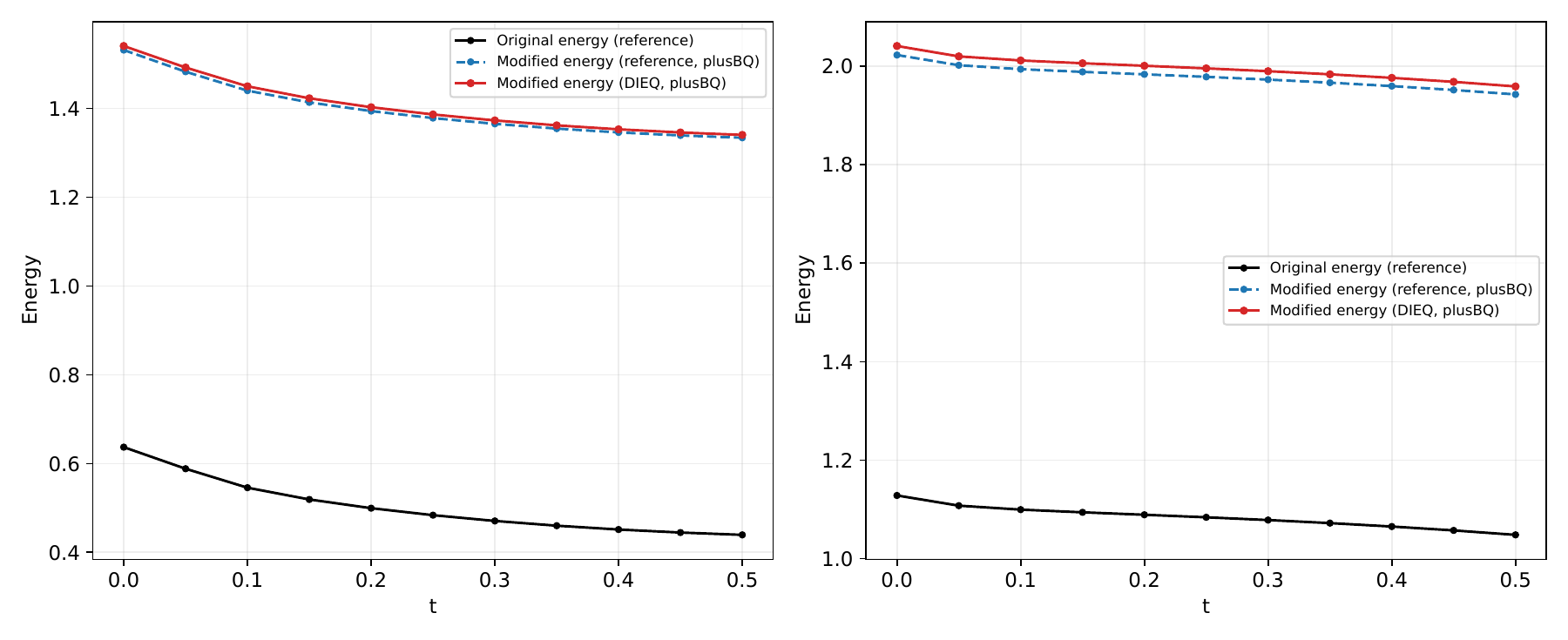}
    \caption{Comparison of the original and modified reference energies with the predicted modified energy for the two-droplet problems. Left: coarsening; right: coalescence.}
    \label{fig:cacns-energy-evolution}
\end{figure}

For the coarsening problem, the smaller droplet gradually shrinks while the larger droplet grows, as illustrated by the phase-field snapshots in Fig.~\ref{fig:cacns-coarsening}. For the coalescence problem, the two initially separated droplets merge into a single droplet and gradually evolve toward a more stable configuration, as shown in Fig.~\ref{fig:cacns-coalescence}. In both cases, the predicted solutions exhibit good agreement with the reference solutions, capturing the characteristic interfacial dynamics while preserving the total phase mass. Fig.~\ref{fig:cacns-energy-evolution} compares the modified energy of the recursive network predictions, evaluated at intervals of $dt=0.05$, with the original and modified reference energies. For both coarsening and coalescence, the predicted modified energy decreases strictly between consecutive recorded time points and remains reasonably close to the modified reference energy. These results provide numerical evidence that the recursive predictions are consistent with the expected energy-dissipation behavior over the time interval $t\in[0,0.5]$.

\section{Conclusion}\label{sec:conclusion}

In this work, we have proposed a class of EDP-PI-DOSnets for solving gradient flow equations. By imposing appropriate constraints on the convolution kernels, we guarantee the symmetric negative-semidefinite structure of the learned operators. Classical energy-stable schemes, including CN, SAV/CN, and R-DVD-1, are then incorporated into the network architecture, ensuring the corresponding discrete energy dissipation properties without introducing additional energy-penalty terms in the loss function. Furthermore, a learnable extrapolation strategy is introduced, and the structure-preserving framework is extended to the CAC-NS system. Numerical experiments demonstrate that the proposed method achieves satisfactory prediction accuracy while preserving the corresponding discrete energy dissipation laws and maintaining good stability during long-time prediction. These results suggest that incorporating classical structure-preserving numerical schemes into neural networks provides an effective approach to constructing energy-stable neural solvers.

\section*{Acknowledgment}

\bibliographystyle{unsrt}
\bibliography{Reference}

@book{leveque2007finite,
  title={Finite difference methods for ordinary and partial differential equations: steady-state and time-dependent problems},
  author={LeVeque, Randall J},
  year={2007},
  publisher={SIAM}
}

@book{zienkiewicz2013finite,
  title={The finite element method: its basis and fundamentals},
  author={Zienkiewicz, Olgierd Cecil and Taylor, Robert Leroy and Zhu, Jian Z},
  year={2013},
  publisher={Elsevier}
}

@book{leveque2002finite,
  title={Finite volume methods for hyperbolic problems},
  author={LeVeque, Randall J and others},
  volume={31},
  year={2002},
  publisher={Cambridge University Press}
}

@book{canuto2006spectral,
  title={Spectral methods},
  author={Canuto, Claudio and Hussaini, M Youssuff and Quarteroni, Alfio and Zang, Thomas A},
  volume={285},
  year={2006},
  publisher={Springer}
}

@article{chorin1967numerical,
  title={The numerical solution of the {Navier--Stokes} equations for an incompressible fluid},
  author={Chorin, Alexandre Joel},
  journal={Bulletin of the American Mathematical Society},
  volume={73},
  number={6},
  pages={928--931},
  year={1967}
}

@incollection{patankar1983calculation,
  title={A calculation procedure for heat, mass and momentum transfer in three-dimensional parabolic flows},
  author={Patankar, Suhas V and Spalding, D Brian},
  booktitle={Numerical prediction of flow, heat transfer, turbulence and combustion},
  pages={54--73},
  year={1983},
  publisher={Elsevier}
}

@article{drew1982mathematical,
  title={Mathematical modeling of two-phase flow},
  author={Drew, Donald A},
  journal={Technical Summary Report Wisconsin Univ},
  year={1982}
}

@article{raissi2019physics,
  title={Physics-informed neural networks: A deep learning framework for solving forward and inverse problems involving nonlinear partial differential equations},
  author={Raissi, Maziar and Perdikaris, Paris and Karniadakis, George E},
  journal={Journal of Computational Physics},
  volume={378},
  pages={686--707},
  year={2019},
  publisher={Elsevier}
}

@article{sirignano2018dgm,
  title={{DGM}: A deep learning algorithm for solving partial differential equations},
  author={Sirignano, Justin and Spiliopoulos, Konstantinos},
  journal={Journal of Computational Physics},
  volume={375},
  pages={1339--1364},
  year={2018},
  publisher={Elsevier}
}

@article{yu2018deep,
  title={The deep {Ritz} method: a deep learning-based numerical algorithm for solving variational problems},
  author={E, Weinan and Yu, Bing},
  journal={Communications in Mathematics and Statistics},
  volume={6},
  number={1},
  pages={1--12},
  year={2018},
  publisher={Springer}
}

@article{zang2020weak,
  title={Weak adversarial networks for high-dimensional partial differential equations},
  author={Zang, Yaohua and Bao, Gang and Ye, Xiaojing and Zhou, Haomin},
  journal={Journal of Computational Physics},
  volume={411},
  pages={109409},
  year={2020},
  publisher={Elsevier}
}

@article{lu2021learning,
  title={Learning nonlinear operators via {DeepONet} based on the universal approximation theorem of operators},
  author={Lu, Lu and Jin, Pengzhan and Pang, Guofei and Zhang, Zhongqiang and Karniadakis, George Em},
  journal={Nature Machine Intelligence},
  volume={3},
  number={3},
  pages={218--229},
  year={2021},
  publisher={Nature Publishing Group UK London}
}

@article{li2020fourier,
  title={Fourier neural operator for parametric partial differential equations},
  author={Li, Zongyi and Kovachki, Nikola and Azizzadenesheli, Kamyar and Liu, Burigede and Bhattacharya, Kaushik and Stuart, Andrew and Anandkumar, Anima},
  journal={arXiv preprint arXiv:2010.08895},
  year={2020}
}

@article{li2020neural,
  title={Neural operator: Graph kernel network for partial differential equations},
  author={Li, Zongyi and Kovachki, Nikola and Azizzadenesheli, Kamyar and Liu, Burigede and Bhattacharya, Kaushik and Stuart, Andrew and Anandkumar, Anima},
  journal={arXiv preprint arXiv:2003.03485},
  year={2020}
}

@article{lan2023dosnet,
  title={{DOSnet} as a non-black-box {PDE} solver: When deep learning meets operator splitting},
  author={Lan, Yuan and Li, Zhen and Sun, Jie and Xiang, Yang},
  journal={Journal of Computational Physics},
  volume={491},
  pages={112343},
  year={2023},
  publisher={Elsevier}
}

@article{tripura2023wavelet,
  title={Wavelet neural operator for solving parametric partial differential equations in computational mechanics problems},
  author={Tripura, Tapas and Chakraborty, Souvik},
  journal={Computer Methods in Applied Mechanics and Engineering},
  volume={404},
  pages={115783},
  year={2023},
  publisher={Elsevier}
}

@article{chen2024positional,
  title={Positional knowledge is all you need: Position-induced transformer {(PiT)} for operator learning},
  author={Chen, Junfeng and Wu, Kailiang},
  journal={arXiv preprint arXiv:2405.09285},
  year={2024}
}

@inproceedings{hao2023gnot,
  title={{GNOT}: A general neural operator transformer for operator learning},
  author={Hao, Zhongkai and Wang, Zhengyi and Su, Hang and Ying, Chengyang and Dong, Yinpeng and Liu, Songming and Cheng, Ze and Song, Jian and Zhu, Jun},
  booktitle={International Conference on Machine Learning},
  pages={12556--12569},
  year={2023},
  organization={PMLR}
}

@article{ye2024pdeformer,
  title={{PDEformer}: Towards a foundation model for one-dimensional partial differential equations},
  author={Ye, Zhanhong and Huang, Xiang and Chen, Leheng and Liu, Hongsheng and Wang, Zidong and Dong, Bin},
  journal={arXiv preprint arXiv:2402.12652},
  year={2024}
}

@article{wu2024transolver,
  title={{Transolver}: A fast transformer solver for {PDEs} on general geometries},
  author={Wu, Haixu and Luo, Huakun and Wang, Haowen and Wang, Jianmin and Long, Mingsheng},
  journal={arXiv preprint arXiv:2402.02366},
  year={2024}
}

@article{wang2021learning,
  title={Learning the solution operator of parametric partial differential equations with physics-informed {DeepONets}},
  author={Wang, Sifan and Wang, Hanwen and Perdikaris, Paris},
  journal={Science Advances},
  volume={7},
  number={40},
  pages={eabi8605},
  year={2021},
  publisher={American Association for the Advancement of Science}
}

@article{li2024physics,
  title={Physics-informed neural operator for learning partial differential equations},
  author={Li, Zongyi and Zheng, Hongkai and Kovachki, Nikola and Jin, David and Chen, Haoxuan and Liu, Burigede and Azizzadenesheli, Kamyar and Anandkumar, Anima},
  journal={ACM/IMS Journal of Data Science},
  volume={1},
  number={3},
  pages={1--27},
  year={2024},
  publisher={ACM New York, NY}
}

@article{gao2021phygeonet,
  title={{PhyGeoNet}: Physics-informed geometry-adaptive convolutional neural networks for solving parameterized steady-state {PDEs} on irregular domain},
  author={Gao, Han and Sun, Luning and Wang, Jian-Xun},
  journal={Journal of Computational Physics},
  volume={428},
  pages={110079},
  year={2021},
  publisher={Elsevier}
}

@article{navaneeth2024physics,
  title={Physics informed {WNO}},
  author={Navaneeth, N and Tripura, Tapas and Chakraborty, Souvik},
  journal={Computer Methods in Applied Mechanics and Engineering},
  volume={418},
  pages={116546},
  year={2024},
  publisher={Elsevier}
}

@article{sanz1992symplectic,
  title={Symplectic integrators for {Hamiltonian} problems: an overview},
  author={Sanz-Serna, Jesus M},
  journal={Acta Numerica},
  volume={1},
  pages={243--286},
  year={1992},
  publisher={Cambridge University Press}
}

@article{bridges2001multi,
  title={Multi-symplectic integrators: numerical schemes for {Hamiltonian} {PDEs} that conserve symplecticity},
  author={Bridges, Thomas J and Reich, Sebastian},
  journal={Physics Letters A},
  volume={284},
  number={4--5},
  pages={184--193},
  year={2001},
  publisher={Elsevier}
}

@article{arnold2006finite,
  title={Finite element exterior calculus, homological techniques, and applications},
  author={Arnold, Douglas N and Falk, Richard S and Winther, Ragnar},
  journal={Acta Numerica},
  volume={15},
  pages={1--155},
  year={2006},
  publisher={Cambridge University Press}
}

@article{hyman2002mimetic,
  title={Mimetic finite difference methods for diffusion equations},
  author={Hyman, J and Morel, J and Shashkov, Mikhail and Steinberg, Stanly},
  journal={Computational Geosciences},
  volume={6},
  number={3},
  pages={333--352},
  year={2002},
  publisher={Springer}
}

@article{dahlby2011general,
  title={A general framework for deriving integral preserving numerical methods for {PDEs}},
  author={Dahlby, Morten and Owren, Brynjulf},
  journal={SIAM Journal on Scientific Computing},
  volume={33},
  number={5},
  pages={2318--2340},
  year={2011},
  publisher={SIAM}
}

@article{elliott1993global,
  title={The global dynamics of discrete semilinear parabolic equations},
  author={Elliott, Charles M and Stuart, Andrew M},
  journal={SIAM Journal on Numerical Analysis},
  volume={30},
  number={6},
  pages={1622--1663},
  year={1993},
  publisher={SIAM}
}

@article{eyre1998unconditionally,
  title={Unconditionally gradient stable time marching the {Cahn-Hilliard} equation},
  author={Eyre, David J},
  journal={MRS Online Proceedings Library (OPL)},
  volume={529},
  pages={39},
  year={1998},
  publisher={Cambridge University Press}
}

@article{shen2010numerical,
  title={Numerical approximations of {Allen--Cahn} and {Cahn--Hilliard} equations},
  author={Shen, Jie and Yang, Xiaofeng},
  journal={Discrete Contin. Dyn. Syst},
  volume={28},
  number={4},
  pages={1669--1691},
  year={2010}
}

@article{du1991numerical,
  title={Numerical analysis of a continuum model of phase transition},
  author={Du, Qiang and Nicolaides, Roy A},
  journal={SIAM Journal on Numerical Analysis},
  volume={28},
  number={5},
  pages={1310--1322},
  year={1991},
  publisher={SIAM}
}

@article{huang2025high,
  title={High-order energy stable discrete variational derivative schemes for gradient flows},
  author={Huang, Jizu},
  journal={IMA Journal of Numerical Analysis},
  volume={45},
  number={4},
  pages={2329--2360},
  year={2025},
  publisher={Oxford University Press}
}

@article{yang2016linear,
  title={Linear, first and second-order, unconditionally energy stable numerical schemes for the phase field model of homopolymer blends},
  author={Yang, Xiaofeng},
  journal={Journal of Computational Physics},
  volume={327},
  pages={294--316},
  year={2016},
  publisher={Elsevier}
}

@article{shen2018scalar,
  title={The scalar auxiliary variable {(SAV)} approach for gradient flows},
  author={Shen, Jie and Xu, Jie and Yang, Jiang},
  journal={Journal of Computational Physics},
  volume={353},
  pages={407--416},
  year={2018},
  publisher={Elsevier}
}

@article{li2020sav,
  title={On a {SAV-MAC} scheme for the {Cahn--Hilliard--Navier--Stokes} phase-field model and its error analysis for the corresponding {Cahn--Hilliard--Stokes} case},
  author={Li, Xiaoli and Shen, Jie},
  journal={Mathematical Models and Methods in Applied Sciences},
  volume={30},
  number={12},
  pages={2263--2297},
  year={2020},
  publisher={World Scientific}
}

@article{yang2021novel,
  title={A novel fully decoupled scheme with second-order time accuracy and unconditional energy stability for the {Navier-Stokes} equations coupled with mass-conserved {Allen-Cahn} phase-field model of two-phase incompressible flow},
  author={Yang, Xiaofeng},
  journal={International Journal for Numerical Methods in Engineering},
  volume={122},
  number={5},
  pages={1283--1306},
  year={2021},
  publisher={Wiley Online Library}
}

@article{greydanus2019hamiltonian,
  title={Hamiltonian neural networks},
  author={Greydanus, Samuel and Dzamba, Misko and Yosinski, Jason},
  journal={Advances in Neural Information Processing Systems},
  volume={32},
  year={2019}
}

@article{chen2019symplectic,
  title={Symplectic recurrent neural networks},
  author={Chen, Zhengdao and Zhang, Jianyu and Arjovsky, Martin and Bottou, L{\'e}on},
  journal={arXiv preprint arXiv:1909.13334},
  year={2019}
}

@article{jin2020sympnets,
  title={SympNets: Intrinsic structure-preserving symplectic networks for identifying {Hamiltonian} systems},
  author={Jin, Pengzhan and Zhang, Zhen and Zhu, Aiqing and Tang, Yifa and Karniadakis, George Em},
  journal={Neural Networks},
  volume={132},
  pages={166--179},
  year={2020},
  publisher={Elsevier}
}

@article{matsubara2022finde,
  title={{FINDE}: Neural differential equations for finding and preserving invariant quantities},
  author={Matsubara, Takashi and Yaguchi, Takaharu},
  journal={arXiv preprint arXiv:2210.00272},
  year={2022}
}

@article{cardoso2025exactly,
  title={Exactly conservative physics-informed neural networks and deep operator networks for dynamical systems},
  author={Cardoso-Bihlo, Elsa and Bihlo, Alex},
  journal={Neural Networks},
  volume={181},
  pages={106826},
  year={2025},
  publisher={Elsevier}
}

@article{matsubara2020deep,
  title={Deep energy-based modeling of discrete-time physics},
  author={Matsubara, Takashi and Ishikawa, Ai and Yaguchi, Takaharu},
  journal={Advances in Neural Information Processing Systems},
  volume={33},
  pages={13100--13111},
  year={2020}
}

@article{wu2023energy,
  title={Energy stable neural network for gradient flow equations},
  author={Wu, Yue and Jin, Tianyu and Chen, Chuqi and Fan, Ganghua and Lan, Yuan and Zhang, Luchan and Xiang, Yang},
  journal={arXiv preprint arXiv:2309.10002},
  year={2023}
}

@article{kutuk2025energy,
  title={Energy dissipation preserving physics informed neural network for {Allen--Cahn} equations},
  author={K{\"u}t{\"u}k, Mustafa and Y{\"u}cel, Hamdullah},
  journal={Journal of Computational Science},
  volume={87},
  pages={102577},
  year={2025},
  publisher={Elsevier}
}

@article{luo2026sav,
  title={A {SAV} neural network method for {Allen--Cahn} and {Cahn--Hilliard} equations},
  author={Luo, Tao and Wang, Guihong and Zhao, Wei},
  journal={arXiv preprint arXiv:2610.08088},
  year={2026}
}

@article{hu2024energetic,
  title={Energetic variational neural network discretizations of gradient flows},
  author={Hu, Ziqing and Liu, Chun and Wang, Yiwei and Xu, Zhiliang},
  journal={SIAM Journal on Scientific Computing},
  volume={46},
  number={4},
  pages={A2528--A2556},
  year={2024},
  publisher={SIAM}
}

@article{zhang2024energy,
  title={Energy-dissipative evolutionary deep operator neural networks},
  author={Zhang, Jiahao and Zhang, Shiheng and Shen, Jie and Lin, Guang},
  journal={Journal of Computational Physics},
  volume={498},
  pages={112638},
  year={2024},
  publisher={Elsevier}
}

@article{li2026structure,
  title={{Structure-preserving Fourier Neural Operators for Long-time Cahn-Hilliard Dynamics under Coarse Temporal Supervision}},
  author={Li, Yingli and Zhou, YC},
  journal={arXiv preprint arXiv:2609.39053},
  year={2026}
}

@article{huang2026pidosnet,
  title={{PI-DOSnet}: A physics-informed deep operator-splitting network for evolution partial differential equations},
  author={Huang, Jizu and Qian, Yue and Zhou, Tao},
  journal={arXiv preprint arXiv:2606.22514},
  year={2026}
}

@book{horn2012matrix,
  title={Matrix analysis},
  author={Horn, Roger A and Johnson, Charles R},
  year={2012},
  publisher={Cambridge University Press}
}

@article{gray2006toeplitz,
  title={Toeplitz and circulant matrices: A review},
  author={Gray, Robert M},
  journal={Foundations and Trends in Communications and Information Theory},
  volume={2},
  number={3},
  pages={155--239},
  year={2006},
  publisher={Emerald Publishing Limited}
}

@misc{driscoll2014chebfun,
  title={Chebfun guide},
  author={Driscoll, Tobin A and Hale, Nicholas and Trefethen, Lloyd N and others},
  year={2014},
  howpublished={Pafnuty Publications, Oxford}
}

\end{document}